\documentclass{article}

\usepackage{PRIMEarxiv}

\usepackage[utf8]{inputenc} 
\usepackage[T1]{fontenc}    
\usepackage{hyperref}       
\usepackage{url}            
\usepackage{booktabs}       
\usepackage{amsfonts}       
\usepackage{nicefrac}       
\usepackage{microtype}      
\usepackage{lipsum}
\usepackage{fancyhdr}       
\usepackage{graphicx}       
\graphicspath{{media/}}     
\usepackage{framed,multirow}
\usepackage{amsmath,amsfonts}
\usepackage{algorithmic}
\usepackage{graphicx}
\usepackage{textcomp}
\usepackage[utf8]{inputenc}
\usepackage[english]{babel}
\usepackage{tikz}
\usepackage{dsfont}
\usepackage{amsthm}
\usepackage{enumitem}
\usepackage{mathtools}
\usepackage{eqnarray}
\usepackage{pifont}
\usepackage{subcaption}
\usepackage{titlesec}
\usetikzlibrary{calc}

\newtheorem{remark}{Remark}
\newtheorem{definition}{Definition}[section]

\DeclareMathOperator{\comp}{C}
\DeclareMathOperator{\g}{g}
\DeclareMathOperator{\shockfilter}{sf}
\newcommand{\alphat}{{\alpha_t(i)}} 
\newcommand{\bs}{{b_{s}}}
\newcommand{\bt}{{b_{t}}}
\newcommand{\bts}{{b_{t+s}}}
\newcommand{\btsd}{{b_{t}}}
\newcommand{\bssd}{{b_{s}}}
\newcommand{\btstwo}{{b_{t}}}
\newcommand{\btssd}{{b_{t+s}}}
\newcommand{\Rmat}{R_{\mu,x}}

\newcommand{\cD}{\mathcal{D}}
\newcommand{\cF}{\mathcal{F}_\mu}

\newcommand{\gc}{\mathcal{GC}}

\newcommand{\gcmux}{\gc_{\mu,x}}
\renewcommand{\S}{\mathcal{S}^2}

\newcommand{\Iprime}{{I}}

\newcommand{\Sd}{{\mathcal{S}^{d-1}}}
\newcommand{\Id}{I} 

\newcommand{\Rq}{\mathbb{R}^q}
\newcommand{\Rd}{\mathbb{R}^d}
\newcommand{\Zq}{\mathbb{Z}^q}

\newcommand{\Rdd}{\mathbb{R}^{d\times d}}

\newcommand{\N}{\mathbb{N}}
\newcommand{\R}{\mathbb{R}}

\title{Mathematical morphology on directional data
}

\author{ Konstantin Hauch \\
	Department of Mathematics,\\ 
    RPTU Kaiserslautern-Landau,\\ 
    Kaiserslautern, Germany\\
    \texttt{hauch@mathematik.uni-kl.de} \\
	\And
	Claudia Redenbach \\
	Department of Mathematics,\\ 
    RPTU Kaiserslautern-Landau,\\ 
    Kaiserslautern, Germany\\
    \texttt{redenbach@mathematik.uni-kl.de} \\
}

\begin{document}
\maketitle

\begin{abstract}
We define morphological operators and filters for directional images whose pixel values are unit vectors.
This requires an ordering relation for unit vectors which is obtained by using depth functions.
They provide a centre-outward ordering with respect to a specified centre vector. 
We apply our operators on synthetic directional images and compare them with classical morphological operators for grey-scale images. 
As application examples, we enhance the fault region in a compressed glass foam and segment misaligned fibre regions of glass fibre reinforced polymers.
\end{abstract}

\keywords{depth function \and $\mu$CT imaging \and image processing \and filtering \and glass foam \and glass fibre reinforced polymer}

\section{Introduction}
Mathematical morphology combines nonlinear image processing and filtering techniques based on ideas from (random) set theory, topology, and stochastic geometry.
It is widely applied for the analysis of spatial structures, e.g. from geology, biology or materials science. 
The foundations of mathematical morphology are laid in the books by Matheron \cite{Matheron74}, Serra \cite{Serra83,Serra88}, and Soille \cite{Soille03book}.
Sternberg  \cite{Sternberg86} generalised mathematical morphology to numerical functions via the umbra method, i.e., the application of set morphology to the graph of the function.
Ronse \cite{Ronse90} and Goutsias et al. \cite{Goutsias95} derived a further generalisation to complete lattices, that is, partially ordered sets with the property that all subsets have a supremum and an infimum.

The progress of imaging methods from binary via grey-scale images to vector-valued images, such as colour or ultra-spectral images, requires modified morphological operators.
In particular, diffusion tensor imaging (DTI) in neuroimaging \cite{ZHANG2020-Overview} results in vector-valued images showing the diffusion of water molecules in the brain. These images can be used to better understand neural diseases as the diffusion in diseased brains is disturbed or altered compared to healthy brains. For instance, Zhang et al. \cite{ZHANG2021-Alzheimer} proposed a novel Alzheimer’s disease multi-class classification framework with embedding feature selection and fusion based on multi-modal neuroimaging.

In mathematical morphology, the step from univariate to multivariate pixel values involves the challenging task of ordering vectors to define the notions of minimum and maximum of a set of vectors.
An ordering structure for vectors can be derived from the concept of depth \cite{Serfling00}.
For a set of vectors, depth functions assign to each vector a value that measures its "centrality" within the set.
The "centre" is the vector maximising the depth function.
Thus, a centre-outward ordering based on the depth values of each vector yields a sound definition of minimum and maximum.

Velasco-Forero et al. \cite{Angulo12depth} defined multivariate mathematical morphological operators via random projection depth. 
They illustrated their approach on colour and hyperspectral images.
Concepts of depth for directional data, i.e. unit vectors in $\Rd$, have been studied, for instance, by Liu et al. \cite{Liu92} and Pandolfo et al. \cite{Pandolfo17}.
Ley et al. \cite{Ley14} introduced quantiles for directional data and the angular Mahalanobis depth, and Garc\'{i}a-Portugu\'{e}s et al. \cite{Verdebout20} developed optimal tests for rotational symmetry against new classes of hyperspherical distributions.

Based on an ordering on the unit sphere $\Sd$, mathematical morphology was extended to directional images by several authors.
Roerdink \cite{Roerdink90} introduced mathematical morphology on the sphere via generalised Minkowski operations.
His motivation is morphological processing of earth pictures under consideration of the curvature. 
Morphological operators for angle-valued images were introduced by Peters \cite{Peters97} and Hanbury et al. \cite{Hanbury01}.
However, a generalisation of their results to the unit sphere $\Sd$ with $d> 2$ is not straightforward.
Frontera-Pons and Angulo \cite{Angulo12morphsphere} defined morphological operators via a local partial ordering. 
They used the Fr\'{e}chet-Karcher barycenter as local origin $\mu$ on the sphere. 
The maximum and minimum of a set of vectors are then found via projecting the vectors into the tangent space at $\mu$.
A drawback of their approach is that the minimum and maximum derived this way are not necessarily elements of the given vector set, which seems unnatural.
Angulo \cite{Angulo15scale-space} extended the concept of multi-scale morphological operators, i.e. a scale-dependent morphology and filters, from $\R$-valued to $\Rd$-valued images. He defined morphological scale-space operators on metric Maslov-measurable spaces for images supported on point clouds. An application of his operators to RGB-valued point clouds, e.g. for object extraction, showed the usefulness of applying the theory of multi-scale morphological operators to $\Rd$-valued images.

In this work, we introduce morphological operators on directional images. 
Section \ref{sec:Basics and Notation} contains basics about mathematical morphology for grey-scale images and statistical depth functions on $\Sd$.
In Section \ref{sec:Mathematical morphology on directional data using directional projection depth}, we introduce the concept of mathematical morphology on directional data using the directional projection depth.
These morphological operators are extended to multi-scale operators in Section \ref{sec:Morphological multi-scale operators for directional images}.
We illustrate and interpret our findings on synthetic and real-world $\S$-valued images in Section \ref{sec:Examples}. A short conclusion is given in Section \ref{sec:Conclusion}.

\section{Basics and Notation}\label{sec:Basics and Notation}
\subsection{Images}
We denote by $\Sd = \{x\in \Rd : \lVert x \rVert _2 = 1\}$ the $(d-1)$-dimensional unit sphere and $\lVert x \rVert _2 = \sqrt{x^T x}$ the Euclidean norm.
An image is a mapping $I:E \rightarrow S$ that is defined on a $q$-dimensional domain $E$ with $E\subset \Rq$ or $E\subset \Zq$. The set $S$ is the set of pixel values. We call $I$ a binary image if $S=\{0,1\}$, a grey-scale image if $S\subset\R$, a vector-valued image if $S\subset\Rd$ and a directional ($\Sd$-valued) image if $S \subset \Sd$.
Usually, $q,d\in\{2,3\}$.
In this paper, we consider examples of $q$-dimensional $\Sd$-valued images with $q=2,3$ and $d=3$. 
We denote by 
\begin{itemize}
    \item $\max I$ the maximal possible pixel value of an image $I$. For instance, $\max I=1$ if $I$ is a binary image or $\max I=255$ if $I$ is an 8-bit grey-scale image.
    \item $\overline\R = \R \cup \{-\infty, \infty\}$ the extended real numbers.
    \item $SO(d)$ the rotation group on $\Rd$.
     \item $x\times y$ the cross product between two vectors $x,y\in\Rd$.
     \item $\gcmux = \{y\in\Rd : (x\times\mu)^Ty=0\} \cap \Sd$
    the great circle containing the vectors $\mu$ and $x\in\Sd$.
    For $d=3$, the great circle $\gcmux$ corresponds to a closed curve on the surface of $\S$ created by the intersection of $\S$ and a $2$-dimensional hyperplane $H$ passing through the origin $0_3$ \cite{FLE87}.
    \item $\check{b}$ the reflection of a function $b:\Rq\rightarrow \overline{\R}$, i.e., $\check{b}(i) = b(-i)$ for all $i\in \Rq$, and by $\check{B}$ the reflection of a set $B\subset\Rq$, i.e., $\check{B} = -B$. A function $b$ is called symmetric if $b=\check{b}$ and a set $B$ is symmetric if $B=\check{B}$.
\end{itemize}

\subsection{Mathematical morphology}
\label{subsec:Morphology}
We first define morphological operators on grey-scale images $I$.
See \cite{Matheron74,Serra83,Soille03book} for a detailed introduction.
The two fundamental operations of mathematical morphology are erosion $\Tilde{\varepsilon}_b$ and dilation $\Tilde{\delta}_b$. They depend on a structuring function $b:\Rq\rightarrow\overline\R$ and are defined as
\begin{align}
    \Tilde{\varepsilon}_b(I)(i) &= \inf_{j\in E} \{I(j)-b(j-i)\},  \quad i\in E, 
    \label{eq:erosion}\\
    \Tilde{\delta}_b(I)(i)      &= \sup_{j\in E} \{I(j)+b(j-i)\},  \quad i\in E.
    \label{eq:dilation}
\end{align} 
In the discrete case, supremum and infimum can be replaced by maximum and minimum. We will not discuss edge effects here (see \cite{Jackway97book}).
Throughout the paper, we assume that $b$ is symmetric.

A well-known and widely used example for $b$ is the flat structuring function. Given a structuring element $B\subset \Rq$, it is defined as 
\begin{align*}
    b(i)=\begin{cases}
    0       & i\in B \\
    -\infty & i\not\in B. \\
    \end{cases}
\end{align*}
We assume that $B$ is centred at the origin and symmetric.
Morphological operators with flat structuring functions are called flat operators. Symbols for flat operators will be written with index $B$ rather than $b$, e.g., a flat erosion with structuring element $B$ is denoted by $\Tilde\varepsilon_B$.

Non-flat or volumic structuring functions assign non-constant weights to the pixel values \cite{Serra83}. 
For instance, 
\begin{align*}
    b(i) &= -\left(\frac{\lVert i \rVert _2}{2}\right)^2
\end{align*}
is a non-flat structuring function.
Note that the grey-scale ranges of images dilated or eroded by non-flat structuring functions are not bounded \cite{Soille03book}.
For instance, an erosion with a non-flat structuring function of an image with non-negative pixel values can result in negative pixel values.
Flat structuring functions do not suffer from that problem.
The output of a flat erosion or dilation is bounded by the grey-scale range of the input image.

The composition of erosion and dilation yields the morphological operators opening $\Tilde{\gamma}_b$ and closing $\Tilde{\varphi}_b$
\begin{align*}
    \Tilde{\gamma}_b(I)(i) &= \Tilde{\delta}_b(\Tilde{\varepsilon}_b(I)) (i),    \\
    \Tilde{\varphi}_b(I)(i)      &= \Tilde{\varepsilon}_b(\Tilde{\delta}_b(I)) (i).
\end{align*} 
In applications, the opening is used to remove bright noise and the closing to remove dark noise.

Some properties of the morphological operators are summarised in the following, see also \cite{Matheron74,Serra83,Jackway97book}. 
Let $i\in E$, $I,I'$ grey-scale images and $\{I_l\}_{l\in\N}$ a family of grey-scale images.
We write $I \leq I'$ if $I(i) \leq I'(i)$ for all $i\in E$. 
Furthermore, $\bigvee$ denotes the point-wise maximum operator and $\bigwedge$ the point-wise minimum operator.
\begin{enumerate}
    \item All the morphological operators are non-linear, i.e., in general
        \begin{align*}
            \Tilde{\Psi}(aI+bI') \not= a\Tilde{\Psi}(I)+b\Tilde{\Psi}(I'),
        \end{align*}
        where $a,b\in\R$ and 
        $\Tilde{\Psi}=\Tilde{\delta}_b, \Tilde{\varepsilon}_b, \Tilde{\gamma}_b$, or $\Tilde{\varphi}_b$.
    \item Dilation and erosion are dual w.r.t. complementation $\Tilde\comp$ with $\Tilde\comp I(i)=\max I -I(i)$, i.e.,
        \begin{align*}
            \Tilde{\delta}_b (I) 
            &= \Tilde\comp \Tilde{\varepsilon}_b( \Tilde\comp I)  
        \end{align*}
        with $b(i) = b(-i)$ due to symmetry.
    \item Opening and closing are dual w.r.t. complementation $\Tilde\comp$, i.e.
        \begin{align*}
            \Tilde{\gamma}_b (I) 
            &= \Tilde\comp \Tilde{\varphi}_b( \Tilde\comp I) 
        \end{align*}
    \item Opening and closing are idempotent, i.e., 
        \begin{align*}
            \Tilde{\gamma}_b (\Tilde{\gamma}_b (I)) &= \Tilde{\gamma}_b (I), \\ 
            \Tilde{\varphi}_b (\Tilde{\varphi}_b (I)) &= \Tilde{\varphi}_b (I).
        \end{align*}
    \item The following distribution laws hold
        \begin{align*}
            \Tilde{\delta}_b \left(\bigvee_l I_l\right) 
            &= \bigvee_l \Tilde{\delta}_b (I_l), \\
            \Tilde{\varepsilon}_b \left(\bigwedge_l I_l\right) 
            &= \bigwedge_l  \Tilde{\varepsilon}_b (I_l).
        \end{align*}
    \item Dilation is associative and erosion fulfils the chain rule, i.e.,
        \begin{align*}
            \Tilde{\delta}_b (\Tilde{\delta}_{b'} (I))  
            &= \Tilde{\delta}_{\Tilde{\delta}_b(b')} (I) ,   \\
           \Tilde{\varepsilon}_b (\Tilde{\varepsilon}_{b'} (I)) 
            &= \Tilde{\varepsilon}_{\Tilde{\delta}_b(b')} (I) .
        \end{align*}
    \item All morphological operators are increasing, i.e.,
        \begin{align*}
            I \leq I' \Rightarrow \Tilde{\Psi}(I) \leq \Tilde{\Psi}(I') ,
        \end{align*}
        where 
        $\Tilde{\Psi}=\Tilde{\delta}_b, \Tilde{\varepsilon}_b, \Tilde{\gamma}_b$, or $\Tilde{\varphi}_b$.
    \item If $b$ is defined at the origin and $b(0)\geq 0$, dilation is extensive and  erosion is anti-extensive, i.e.,
        \begin{equation*}
            \Tilde{\varepsilon}_b (I) \leq I 
            \leq \Tilde{\delta}_b (I). 
        \end{equation*}
    \item Closing is extensive and opening is anti-extensive, i.e.,
        \begin{equation*}
            \Tilde{\gamma}_b (I)  \leq I \leq \Tilde{\varphi}_b (I).  
        \end{equation*}
\end{enumerate}
The latter two properties can be combined to give the ordering relation
\begin{align}
\label{eq:order}
    \Tilde{\varepsilon}_b(I) &\leq \Tilde{\gamma}_b(I) \leq I \leq \Tilde{\varphi}_b(I) \leq \Tilde{\delta}_b(I).
    \end{align}


\subsection{Morphological scale-space}\label{sec:Morphological scale-space}
In the context of mathematical morphology for grey-scale images, erosion and dilation can be written as scale-space operators \cite{Serra83,Heijmans94}
\begin{align}
    \label{eq:ss-erosion}
	\Tilde{\varepsilon}_\bt(I)(i) &= \inf_{j\in E} \{I(j)-\bt(j-i)\}    \\
    \label{eq:ss-dilation}
	\Tilde{\delta}_\bt(I)(i)      &= \sup_{j\in E} \{I(j)+\bt(j-i)\},
\end{align} 

with 
a scaled structuring function
\begin{align*}
	\bt:\Rq\rightarrow \overline\R, \quad t\geq 0    
\end{align*}
with scale parameter $t$.
Following \cite[p.15]{Jackway97book}, we assume $\bt$ to have the following properties:
\begin{enumerate}
	\item  $\{\bt\}_t$ is a one-parametric family of convex, continuous, symmetric functions and fulfils the semi-group property, i.e.
	\begin{align}
		\bt \dot{+}\bs &=\bts  \qquad t,s \geq 0,
		\label{eq:structuring_function_semi-group}
	\end{align}
	where $\dot{+}$, $+$ are group operations \cite{Heijmans94}.
	For instance, when choosing $b_t$ to be a flat structuring function with structuring element $tB$, the group operations are the Minkowski addition and the standard addition in $\R$:
	\begin{align*}
		tB\oplus sB
		&=(t+s)B  \qquad t,s \geq 0.
	\end{align*}
	\item $\bt$ is non-positive and monotonically decreasing with a global maximum at the origin of value zero, i.e.,
	\begin{align}
		\bt(i)&\leq 0 \qquad ~~\forall ~t \geq 0, i\in \Rq, 
		\label{eq:structuring_function_non-positive}\\
		\bt(i)& \geq \bt(j) \quad \lVert i \rVert _2 < \lVert j \rVert _2, 
		\label{eq:structuring_function_decreasing}\\
		\bt(0)&=0. 
		\label{eq:structuring_function_zero_in_origin}
	\end{align}
\end{enumerate}

An example of a non-flat scaled structuring function is the Poweroid structuring function \cite[p.14, Def. 5]{Jackway97book}
\begin{align}
	\bt(i) &= - t\left(\frac{\lVert i \rVert _2}{t}\right)^a \qquad a\geq0, t>0.
	\label{eq:Poweroid structuring function}
\end{align}
   
\subsection{Directional data}
Directional data are given as a set of vectors on the unit sphere $\Sd$ for some $d\geq 2$.
The vector $x=(x_1,\dots,x_d)^T\in\Sd$ can be represented as a point on the surface of $\Sd$.
For $x$ in Cartesian coordinates, a representation in angular coordinates for arbitrary $d$ is given in \cite{Blumenson60}.

For $d=3$, we use the usual notation for spherical coordinates 
\begin{align}
x=x(\phi,\theta)&=\begin{pmatrix}
\sin{(\theta)}\cos{(\phi)}\\
\sin{(\theta)}\sin{(\phi)}\\
\cos{(\theta)}
\end{pmatrix}
\label{eq:spherical coordinates}
\end{align}
with co-latitude $\theta\in[0,\pi]$ and longitude $\phi\in[0,2\pi)$.
A random vector $X\in\S$ can be represented in spherical coordinates by random variables $\Theta$ and $\Phi$ with realisations $\theta$ and $\phi$.

\subsection{Statistical depth functions}
Depth functions are applied in multidimensional non-parametric robust data analysis and establish an ordering relation between vectors.
See \cite{Serfling00} for a detailed introduction.
The literature gives several depth functions, for instance, halfspace depth \cite{Tukey75}, simplicial depth \cite{Liu90}, projection depth \cite{Donoho92}, spatial depth \cite{Vardi00}, or the Mahalanobis depth
\cite{Serfling00}.
For defining a depth for directional data, we initially restrict attention to the class $\cF$ of distributions on $\Sd$ with a bounded density that admits a unique modal direction $\mu$. 
Examples of distributions in $\cF$ are the von Mises–Fisher distribution, the Kent distribution or the Bingham distribution \cite{FLE87}. 
We further assume that $\mu$ coincides with the Fisher spherical median \cite{Fisher85}, that is
\begin{align*}
    \mu= \arg\min_{\gamma\in\Sd}E(\arccos(X^T\gamma)).    
\end{align*}
Ley et al. \cite{Ley14} adapted the definition of a depth function given in \cite[Definition 2.1]{Serfling00} for directional data as follows.
\begin{definition}

\label{def:depth_sphere}
A statistical depth function on $\Sd$ is a bounded, non-negative mapping $\cD_\cdot(\cdot): \cF\times\Sd \rightarrow \R$ satisfying
\begin{enumerate}
    \item \label{def:depth-sphere-prop-1}
    $\cD_{F_{AX}}(AX)=\cD_{F_{X}}(X)$ holds for any random vector $X\in\Sd$ and any rotation matrix $A\in\Rdd$.
    \item \label{def:depth-sphere-prop-2}
    $\cD_F(\mu)=\sup_{y\in \Sd}\cD_F(y)$ holds for any $F\in\cF$.
    \item \label{def:depth-sphere-prop-3}
    For any $F\in\cF$, $\cD_F(x)\leq \cD_F(c(\alpha))$ for the unique geodesic $c:[0,1]\rightarrow\Sd$ with $c(0)=\mu$, $c(1)=x$ and $\alpha\in[0,1]$. 
    \item \label{def:depth-sphere-prop-4}
    $\cD_F(-\mu) = 0$ for each $F\in\cF$ where $-\mu$ is the antipodal point of the centre $\mu$.
\end{enumerate}
\end{definition}
Property \ref{def:depth-sphere-prop-1} implies that the depth of a vector should be independent of the underlying coordinate system.
Property \ref{def:depth-sphere-prop-2} means that for any $F$ with unique centre $\mu$, the depth function is maximal at $\mu$.
Property \ref{def:depth-sphere-prop-3} indicates that for any point $X$ moving away from the centre $\mu$ along any geodesic starting at $\mu$, the depth of $X$ decreases monotonically.
Property \ref{def:depth-sphere-prop-4} implies that the depth value of $X$ is zero at $-\mu$ since $-\mu$ is at maximal distance from $\mu$ on $\Sd$.

Depth functions for directional data are, for instance, directional distance-based depths \cite{Pandolfo17} or the angular Mahalanobis depth \cite{Ley14}. 
For defining morphological operators, we cannot use the latter since its application violates the ordering property \eqref{eq:order} for morphological operators. 
We rather use a re-scaled version of the cosine distance depth from \cite{Pandolfo17} as defined in the following.
\subsection{The angular projection depth}
Let $X,X_1,\dots,X_n\in \Sd$ be i.i.d. random vectors with $X\sim F\in\cF$. 
We define a depth for directional data by 
\begin{align}
    D^{proj}_{F}(X) = \frac{1 + X^T\mu}{2}, \quad X\in \Sd.
    \label{eq:D_proj}
\end{align}
$D^{proj}_F(X)$ provides a centre-outward ordering with $D^{proj}_F(\mu) = 1$, $D^{proj}_F(-\mu) = 0$ and is decreasing on a geodesic from $\mu$ to $-\mu$.
Let $A\in\Rdd$ be any rotation matrix.
Then the distribution of the transformed vector $AX$ has centre $A\mu$.
It follows that 
\begin{align*}
D^{proj}_{F_{AX}}(AX) &=\frac{1 + (AX)^T A\mu}{2} 
=\frac{1 + X^TA^TA\mu}{2} \\
&=\frac{1 + X^T\mu}{2} 
=D^{proj}_{F_{X}}(X).
\end{align*}
Thus, the properties of a directional depth function given in Definition \ref{def:depth_sphere} are fulfilled.

The empirical angular projection depth reads
\begin{align*}
    D^{proj}(x) &= \frac{1 + x^T\hat{\mu}}{2},
\end{align*}
where $\hat{\mu}$ is the empirical Fisher spherical median \cite{Fisher85}
\begin{align}
    \hat{\mu}= \arg\min_{\gamma\in\S}\sum_{i=1}^N\arccos(X_i^T\gamma).     
    \label{eq:empirical Fisher spherical median}
\end{align}
\section{Mathematical morphology on directional images using directional projection depth}\label{sec:Mathematical morphology on directional data using directional projection depth}
\begin{figure}
    \centering
    \begin{tikzpicture}[
        point/.style = {draw, circle, fill=black, inner sep=0.7pt},
        ]
        \def\rad{3cm}
        \coordinate (O) at (0,0);  
        \filldraw[ball color=white] (O) circle [radius=\rad];
        \draw[dashed] 
          (\rad,0) arc [start angle=0,end angle=180,x radius=\rad,y radius=5mm];
        \draw
          (\rad,0) arc [start angle=0,end angle=-180,x radius=\rad,y radius=5mm];
        \draw[->, thick](0,0) -- (0,{\rad}) node[above] {$\mu$};
        \draw[-, dashed](0,{-\rad}) -- (0,0);
        \draw[dashed] (O) node[left] {$o$};

        \draw[|-|, thick]({\rad+1cm},-\rad) -- ({\rad+1cm},{\rad}) node[right, yshift = -\rad] {$D^{proj}$};
        \draw ({\rad+1cm},-\rad) node[right] {0};
        \draw ({\rad+1cm},\rad) node[right] {1};
        \def\angletheta_x1{55}
        \draw[dashed] ({\rad*cos(\angletheta_x1)},{\rad*sin(\angletheta_x1)}) node[right] {$x_{3}$};
        \node at ({\rad*cos(\angletheta_x1)},{\rad*sin(\angletheta_x1)}) [circle,fill,inner sep=1.5pt]{};
        \draw[<-, thick](0,{\rad*sin(\angletheta_x1)}) -- ({\rad*cos(\angletheta_x1)},{\rad*sin(\angletheta_x1)});
        \draw[-, thick](0,0) -- ({0},{\rad*sin(\angletheta_x1)}) node[left] {$x_{(3)}^T\mu$};
        \node at ({0},{\rad*sin(\angletheta_x1)}) [circle,fill,inner sep=1.5pt]{};
        \def\angletheta_x2{35}
        \draw[dashed] ({\rad*cos(\angletheta_x2)},{\rad*sin(\angletheta_x2)}) node[right] {$x_{(2)}$};
        \node at ({\rad*cos(\angletheta_x2)},{\rad*sin(\angletheta_x2)}) [circle,fill,inner sep=1.5pt]{};
        \draw[<-, thick](0,{\rad*sin(\angletheta_x2)}) -- ({\rad*cos(\angletheta_x2)},{\rad*sin(\angletheta_x2)});
        \draw[-, thick](0,0) -- ({0},{\rad*sin(\angletheta_x2)}) node[left] {$x_{(2)}^T\mu$};
        \node at ({0},{\rad*sin(\angletheta_x2)}) [circle,fill,inner sep=1.5pt]{};
        \def\angletheta_x3{15}
        \draw[dashed] ({\rad*cos(\angletheta_x3)},{\rad*sin(\angletheta_x3)}) node[right] {$x_{(1)}$};
        \node at ({\rad*cos(\angletheta_x3)},{\rad*sin(\angletheta_x3)}) [circle,fill,inner sep=1.5pt]{};
        \draw[<-, thick](0,{\rad*sin(\angletheta_x3)}) -- ({\rad*cos(\angletheta_x3)},{\rad*sin(\angletheta_x3)});
        \draw[dashed](0,0) -- ({0},{\rad*sin(\angletheta_x3)}) node[left] {$x_{(1)}^T\mu$};
        \node at ({0},{\rad*sin(\angletheta_x3)}) [circle,fill,inner sep=1.5pt]{};
    \end{tikzpicture}
\caption{Vector ordering by directional projection depth.     
}
\label{fig:projection-ordering}
\end{figure}
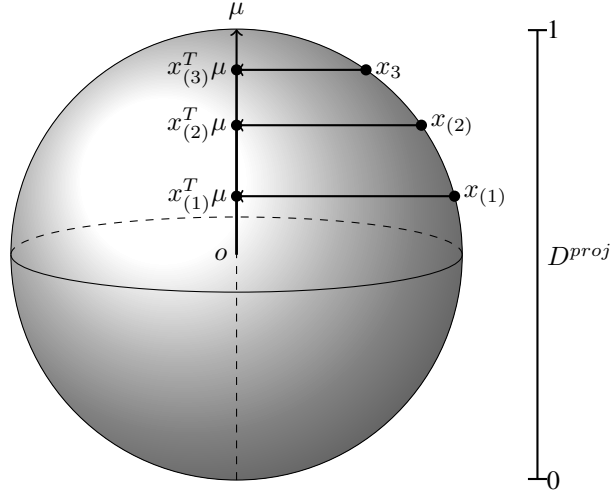    
The definition of the standard operators of mathematical morphology requires an ordering relation between pixel values, see \eqref{eq:erosion} and \eqref{eq:dilation}.
Multivariate extensions of ordering and mathematical morphology were subject to several studies \cite{Angulo07,Aptoula07,Barnett76,Goutsias95,Najman10}. In particular, orderings can be derived from depth functions.

For $\Sd$ valued images, we use the projection depth $D^{proj}$ to define an ordering relation between unit vectors. 
This allows for a definition of erosion, dilation, and the morphological operators derived from them as described below. For simplicity, we concentrate on flat morphological operators with structuring element $B$. 
In a given application, the assumption $I(i)\sim F \in \cF$ may not be fulfilled. 
Nevertheless, a suitable central direction $\mu$ may be derived from the experimental setup.
Hence, $\mu$ will be interpreted as a parameter in the following. We will write $D_\mu^{proj}$ for the projection depth using $\mu \in \Sd$ as the centre.
The parameter $\mu$ can be selected globally for the whole image or in a locally adaptive manner for parts of the image. The latter is not considered here.

The erosion $\varepsilon_B$ of a directional image $\Id$ at pixel position $i\in E$ is (implicitly) defined by
\begin{align*}
    D_\mu^{proj}(\varepsilon_B(\Id)) (i) &= \Tilde{\varepsilon}_B(D_\mu^{proj}(\Id))(i),
\end{align*}
where 
$\Tilde{\varepsilon}$ is the flat erosion of a grey-scale image. 
Analogously, the dilation $\delta_B$ of a directional image $\Id$ at pixel position $i\in E$ is (implicitly) defined by
\begin{align*}
    D_\mu^{proj}(\delta_B(\Id))(i) &= \Tilde{\delta}_B(D_\mu^{proj}(\Id))(i),
\end{align*}
where 
$\Tilde{\delta}$ is the flat dilation of a grey-scale image. 
Hence, a dilation will select the most central direction covered by the structuring element while an erosion selects the most outlying direction.

The explicit definition reads
\begin{align*}
    \varepsilon_B(\Id)(i) &= {D_\mu^{proj}}^{-1}\left(\Tilde{\varepsilon}_B(D_\mu^{proj}(\Id))\right)(i),\\
    \delta_B(\Id)(i) &= {D_\mu^{proj}}^{-1}\left(\Tilde{\delta}_B(D_\mu^{proj}(\Id))\right)(i),
\end{align*}
where ${D_\mu^{proj}}^{-1}$ refers to the preimage under ${D_\mu^{proj}}$.

As $D_\mu^{proj}$ is not injective, the preimage may consist of more than one element. Consider vectors $x_1,x_2\in\Sd$ with $x_1^T\mu=x_2^T\mu$. 
Then, $D_\mu^{proj}(x_1)=D_\mu^{proj}(x_2)$ but $x_1=x_2$ is not necessarily the case.
To resolve this issue we use a lexicographic order for vectors of the same depth value which yields a total ordering \cite{Angulo12depth}.
Rotate the vectors such that $\mu=(0,0,1)^T$.
Among the (pairwise different) vectors $\{x_i\}_{i=1,\dots,n}$ with $D_\mu^{proj}(x_i)=D_\mu^{proj}(x_j)$, we choose the vector $x^*$ with the smallest longitude angle, i.e., 
$x^*=\arg\min\limits_{x_i}\phi_i$ with $(\theta_i,\phi_i)$ spherical coordinates of $x_i$, $i=1,\dots,n$.
Then rotate $x^*$ back.

Opening and closing are defined by
\begin{align*}
    \gamma_B(\Id)(i) &= \delta_B(\varepsilon_B(\Id))(i),\\
    \varphi_B(\Id)(i) &= \varepsilon_B(\delta_B(\Id))(i).
\end{align*}

The theory of $h$-orderings \cite{Goutsias95} immediately implies that the operators defined above fulfil the properties of morphological operators stated in Subsection \ref{subsec:Morphology}, see \cite{Angulo07}. In this theory, $h$ is a map from $\Sd$ into a space with a partial order. Here we choose $h=D^{proj}$.

Of course, further morphological filters have their depth analogy.  
For instance, the scalar difference between the depth of dilation and erosion defines the morphological gradient
\begin{align*}
    \g_{D_\mu^{proj},B}(\Id)(i) 
    &= D_\mu^{proj}\left(\delta_B(\Id)(i)\right) - D_\mu^{proj}\left(\varepsilon_B(\Id)(i)\right).  
\end{align*}
The morphological Laplacian is defined by
\begin{align*}
    \Delta_{D_\mu^{proj},B}(\Id)(i) &= \Delta_\delta(i) - \Delta_\varepsilon(i)
\end{align*}
with 
$\Delta_\delta(i) = D_\mu^{proj}\left(\delta_B(\Id)(i)\right) - D_\mu^{proj}(\Id)(i)$ 
and
$\Delta_\varepsilon(i) = D_\mu^{proj}(\Id)(i) - D_\mu^{proj}\left(\varepsilon_B(\Id)(i)\right)$.
The shock filter is defined by
\begin{align*}
    \shockfilter_{D_\mu^{proj},B}(\Id)(i)  
    =\begin{cases}
    \varepsilon_B(\Id)(i) & \Delta_{D_\mu^{proj},B}(\Id)(i) < 0, \\
    \delta_B(\Id)(i) & \Delta_{D_\mu^{proj},B}(\Id)(i) > 0, \\
    \Id(i) & \text{ otherwise.} \\
    \end{cases}
\end{align*}
Shock filtering is used to enhance edges.
For grey-scale images, the idea is to dilate near local maxima and erode near local minima. 
This way, the contrast between regions of pixel values with high and low depth is enhanced. 


\section{Pseudo morphological multi-scale operators for directional images}\label{sec:Morphological multi-scale operators for directional images}

In the next step, we want to extend morphological operators to a multi-scale setting, leading to a morphological scale space (see \cite{Angulo15scale-space}). A very simple approach to scaling is scaling of a flat structuring element $B$ by a factor $t$. Here, we will consider an alternative approach which yields a structuring function which has a meaningful interpretation, fulfils the semi-group property (at least partially), and results in bounded non-flat operators. We call them pseudo morphological operators since their behaviour reminds of erosion and dilation as demonstrated in Section \ref{sec:Examples}.



\subsection{Structuring function for directional images}

In general, structuring functions $\bt:E\rightarrow \overline\R$ are unbounded.
Thus, multi-scale dilation and erosion inherit this unboundedness if we use $D_\mu^{proj}(\Id(i))+\bt(i)$ and $D_\mu^{proj}(\Id(i))-\bt(i)$, respectively.
However, we cannot associate a vector with a depth value outside $[0,1]$ since no unit vector could be found as preimage.
In our approach we therefore use structuring functions $\bt:E\rightarrow SO(d)$.
We will define the operation $\Id(i)+\bt(i)$ and $\Id(i)-\bt(i)$ to be a vector rotation of $\Id(i)$ where the rotation angle depends on the scale $t$.
Therefore, the depths of $\Id(i)+\bt(i)$ and $\Id(i)-\bt(i)$ are both bounded.

The main idea is to rotate a pixel value $x=\Id(i)$, $i\in E$, about $x\times\mu$ towards $\mu$ or away from $\mu$ on a great circle $\gcmux$.
This increases or decreases the depth value of the rotated $x$ which gives our approach also a meaningful interpretation.
Note that $\Id(i)+\bt(i)$ and $\Id(i)-\bt(i)$, $i\in E$, result again in directional images which seems natural. 

\begin{figure}
	\centering
	\begin{tikzpicture}
		\def\rad{4cm}
		\def\angletheta{135}
		\def\anglethetatowards{110}
		\def\anglethetaaway{160}
		\coordinate (O) at (0,0);  
		
		\begin{scope}
			\clip (-\rad,-1) rectangle (\rad,\rad);
			\draw (0,0) circle(\rad);
		\end{scope}
		
		\draw[->, thick](0,0) -- (0,{\rad}) node[above] {$\mu$};
		\draw[dashed] ({\rad/2^.5},{\rad/2^.5}) node[right] {$\gcmux$};
		\draw[dashed] (0,-0.5) node[below] {$\mathcal{S}^2$};
		
		\draw[->,thick,red] ([shift=(\angletheta:\rad)]0,0) arc (\angletheta:\anglethetatowards:\rad);
		
		\draw[->,thick,red] ([shift=(\angletheta:\rad)]0,0) arc (\angletheta:\anglethetaaway:\rad);
		
		\draw[-,thick,red] ([shift=(\angletheta:{\rad/3})]0,0) arc (\angletheta:\anglethetaaway:{\rad/3})
		node[left, yshift = 0.55cm, xshift = -0.15cm] 
		{$\alpha_{+}$};
		\draw[-,thick,red] ([shift=(\angletheta:{\rad/3})]0,0) arc (\angletheta:\anglethetatowards:{\rad/3})
		node[left, yshift = 0.25cm, xshift = -0.25cm] 
		{$\alpha_{-}$};
		\draw[-,thick] ([shift=(\angletheta:{3cm})]0,0) arc (\angletheta:90:{3cm})
		node[below, yshift = 0cm, xshift = -1.75cm] {$\theta$};
		
		
		\draw[dashed] 
		(\rad,0) arc [start angle=0,end angle=180,x radius=\rad,y radius=5mm];
		\draw
		(\rad,0) arc [start angle=0,end angle=-180,x radius=\rad,y radius=5mm];
		\draw[->, thick](0,0) -- ({\rad*cos(\angletheta)},{\rad*sin(\angletheta)}) node[left] {$x$};
		
		\draw[dashed] 
		(\rad,0) arc [start angle=0,end angle=180,x radius=\rad,y radius=5mm];
		\draw
		(\rad,0) arc [start angle=0,end angle=-180,x radius=\rad,y radius=5mm];
		\draw[->, thick](0,0) -- ({\rad*cos(\anglethetatowards)},{\rad*sin(\anglethetatowards)}) node[above left] 
		{$x-\btsd$};
		
		\draw[dashed] 
		(\rad,0) arc [start angle=0,end angle=180,x radius=\rad,y radius=5mm];
		\draw
		(\rad,0) arc [start angle=0,end angle=-180,x radius=\rad,y radius=5mm];
		\draw[->, thick](0,0) -- ({\rad*cos(\anglethetaaway)},{\rad*sin(\anglethetaaway)}) node[left] 
		{$x+\btsd$};
		
	\end{tikzpicture}
	\caption{
		Illustration of the operations $x - \btsd$ given in \eqref{eq:minus-operation-multi-scale-erosion} and $x + \btsd$ given in \eqref{eq:plus-operation-multi-scale-dilation} with $d=3$ and $x=x(\phi,\theta)$ in spherical coordinates \eqref{eq:spherical coordinates}. 
		$x + \btsd$ increases the geodesic distance between $x$ and $\mu$ by $\alpha_+$.
		$x - \btsd$ decreases the geodesic distance between $x$ and $\mu$ by $\alpha_-$.
		$x$ moves on the great circle $\gcmux$.
	}
	\label{fig:structuring function_construction}
\end{figure}
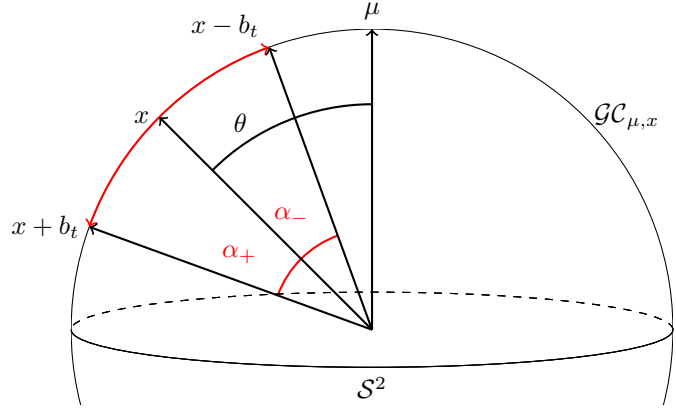
For a formal definition, let $\mu\in\Sd$, $i\in E$, $x=\Id(i)\in\Sd$ with $x^T\mu = \cos(\theta)$ and $\Rmat(\alpha)\in SO(d)$ be a rotation matrix which rotates $x$ about the cross product $\mu\times x$ and $\Rmat(\theta)x=\mu$.
We define a structuring function $\btsd$ for directional images by a mapping
\begin{align}
	\btsd: E   &\rightarrow SO(d) \nonumber\\
	i       &\mapsto    \btsd(i) := \Rmat(\alphat), \label{eq:bt-def}
\end{align}
with rotation angle $\alphat$, $t\geq 0$.

Since $\Rmat$ depends on $x$, the structuring function is locally adaptive.
We choose 
\begin{align}
    \alphat &= \min\left(\frac{\lVert i \rVert _2^2}{t}, \pi \right),
	\label{eq:alpha}
\end{align}
where the truncation at $\pi$ is discussed in more detail in Remark \ref{rem:Fixation of a vector}. Note the similarity of $\alpha_t$ with quadratic scaling as discussed in \cite{Heijmans02}. 


Let $f\Box g$ be the infimal convolution of two functions $f,g: \Rq\to\R$ \cite{Stroemberg1996}, i.e.,
\begin{align}
	(f \Box g) (i) = \inf_{j\in \Rq} \left\{ f(i-j) + g(j) \right\}.
	\label{eq:infimal convolution}
\end{align}
Define
\begin{align}
	(\btsd \dot{+} \bssd)(i) := \Rmat( (\alpha_t \Box \alpha_s)(i)).
	\label{eq:btsd(i) dot+ bssd(i)}
\end{align}	
The structuring element $\btsd$ fulfils the semi-group property  \eqref{eq:structuring_function_semi-group}, i.e. $\btsd(i) \dot{+} \bssd(i) = \btssd (i)$, for $\lVert i \rVert _2^2 / t,\lVert i \rVert _2^2 / s < \pi$ as shown in the Appendix.


We illustrate the construction for $d=3$ and $x=x(\phi,\theta)$ in spherical coordinates as given in \eqref{eq:spherical coordinates}: 
The rotation matrix $\Rmat$ rotates $x$ on the great circle $\gcmux$.
We use $\Rmat$ to move a pixel value $x\in\S$ towards $\mu$ or away from $\mu$ on $\gcmux$.
Thus, the rotation matrix $\Rmat$ causes a change of the angle between $x$ and $\mu$, here the co-latitudinal angle $\theta$. 
The longitude angle $\phi$ is not changed. 
Therefore, we define a rotation of the vector $x$ depending on the structuring function $\btsd$ to this change of $\theta$ as follows (see Figure \ref{fig:structuring function_construction}):
\begin{align}
    x(\phi,\theta) - \btsd(i) &:= x(\phi, \theta-\alpha_-)
	\label{eq:minus-operation-multi-scale-erosion}
\end{align}
with 
\begin{align*}
	\alpha_-&=\begin{cases}
		\alphat,&\text{for } \alphat \in [0,\theta]\\
		\theta, &\text{for } \alphat > \theta. 
	\end{cases}
\end{align*}
\begin{align}
    x(\phi,\theta) + \btsd(i) &:= x(\phi, \theta+\alpha_+)
	\label{eq:plus-operation-multi-scale-dilation}
\end{align}
with  
\begin{align*}
	\alpha_+&=\begin{cases}
		\alphat,   &\text{for } \alphat \in [0,\pi-\theta]\\
		\pi - \theta,    &\text{for } \alphat > \pi-\theta.
	\end{cases}
\end{align*}
$x - \btsd$ has a smaller geodesic distance to $\mu$ than $x$, and $x + \btsd$ has a larger geodesic distance to $\mu$ than $x$.

\begin{remark}[Truncation of rotation]\label{rem:Fixation of a vector}
The restriction of $\alpha_-$ to $\theta $ and  
of $\alpha_+$ to $\pi-\theta $ are necessary to achieve an analogous ordering relation as in \eqref{eq:order}.
Otherwise, $x$ can be rotated beyond $\mu$ which would decrease $D_\mu^{proj}(x - \btsd(i))$.
Analogously, $x$ can be rotated beyond $-\mu$ which would increase $D_\mu^{proj}(x + \btsd(i))$.
To avoid this, we fix the vector at $\mu$ or $-\mu$ respectively, such that no roation beyond is possible. 
\end{remark}

With truncation, the semi-group property is no longer valid.
However, truncation ensures that 
\begin{align}
	D_\mu^{proj}(x + \btsd) &\leq D_\mu^{proj}(x) \leq D_\mu^{proj}(x - \btsd).
	\label{eq:ordering-relation-structuring-element}
\end{align}
We obtain (in some sense) analogue properties to \eqref{eq:structuring_function_non-positive}$-$\eqref{eq:structuring_function_zero_in_origin}:
\begin{enumerate}
	\item $\btsd$ is non-positive w.r.t. $D_\mu^{proj}$ in the sense that for rotations of $x$ away from $\mu$ its depth decreases, i.e.
	\begin{align*}
		D_\mu^{proj}(x + \btsd(i)) - D_\mu^{proj}(x) &\leq 0.
	\end{align*}
	\item $\btsd$ is monotonically decreasing w.r.t. $D_\mu^{proj}$ in the sense that for increasing distance between pixel positions the depth decreases, i.e.
	\begin{align*}
		D_\mu^{proj}(x + \btsd(i))  &\geq D_\mu^{proj}(x + \btsd(j)), \quad \lVert i \rVert _2 < \lVert j \rVert _2.
	\end{align*}
	\item $\btsd$ has a global maximum at the the origin $o\in\Rq$ w.r.t. $D_\mu^{proj}$ with  
	\begin{align*}
		D_\mu^{proj}(x + \btsd(o)) 
  =D_\mu^{proj}(x) 
\end{align*}
\end{enumerate}


\subsection{Pseudo morphological multi-scale operators for directional images}
Let $\btsd$ be as in Equation \eqref{eq:bt-def}.
The multi-scale erosion $\varepsilon_\btsd$ of a directional image $\Id$ at pixel position $i\in E$ is (implicitly) defined by
\begin{align}
	D_\mu^{proj}(\varepsilon_\btsd(\Id)) (i) &= \inf_{j\in E} \{ D_\mu^{proj}(\Id(j) -\btsd(j-i)) \},
	\label{eq:h-erosion-multiscale-implicit}
\end{align}
where $\Id(j) -\btsd(j-i)$ is defined in Equation (\ref{eq:minus-operation-multi-scale-erosion}).
Analogously, the multi-scale dilation $\delta_\btsd$ of a directional image $\Id$ at pixel position $i\in E$ is (implicitly) defined by
\begin{align}
	D_\mu^{proj}(\delta_\btsd(\Id)) (i) &= \sup_{j\in E} \{ D_\mu^{proj}(\Id(j) +\btsd(j-i)) \},
	\label{eq:h-dilation-multiscale-implicit}	
\end{align}
where $\Id(j) +\btsd(j-i)$ is defined in Equation (\ref{eq:plus-operation-multi-scale-dilation}).
The interpretation of the multi-scale operators is as follows:
We want to rotate the vector $\Id(j)$ depending on $||i-j||_2$ and the scale $t$ towards $\mu$ ($D_\mu^{proj}$ increases) or towards $-\mu$ ($D_\mu^{proj}$ decreases). 

Multi-scale opening, closing, morphological gradient and shock filter can be defined as analogues to their flat counterparts.

\section{Examples}\label{sec:Examples}
In this section, we investigate the effect of the morphological operators on synthetic $\S$-valued images and compare them to their standard grey-scale counterparts.
As real application examples, we detect misaligned regions in a glass fibre reinforced composite and 
enhance changes in the displacement field of a compressed glass foam. 
\subsection{Flat morphological operators}
\begin{figure}
    \centering
    \includegraphics[width=0.5\linewidth,angle =-90]{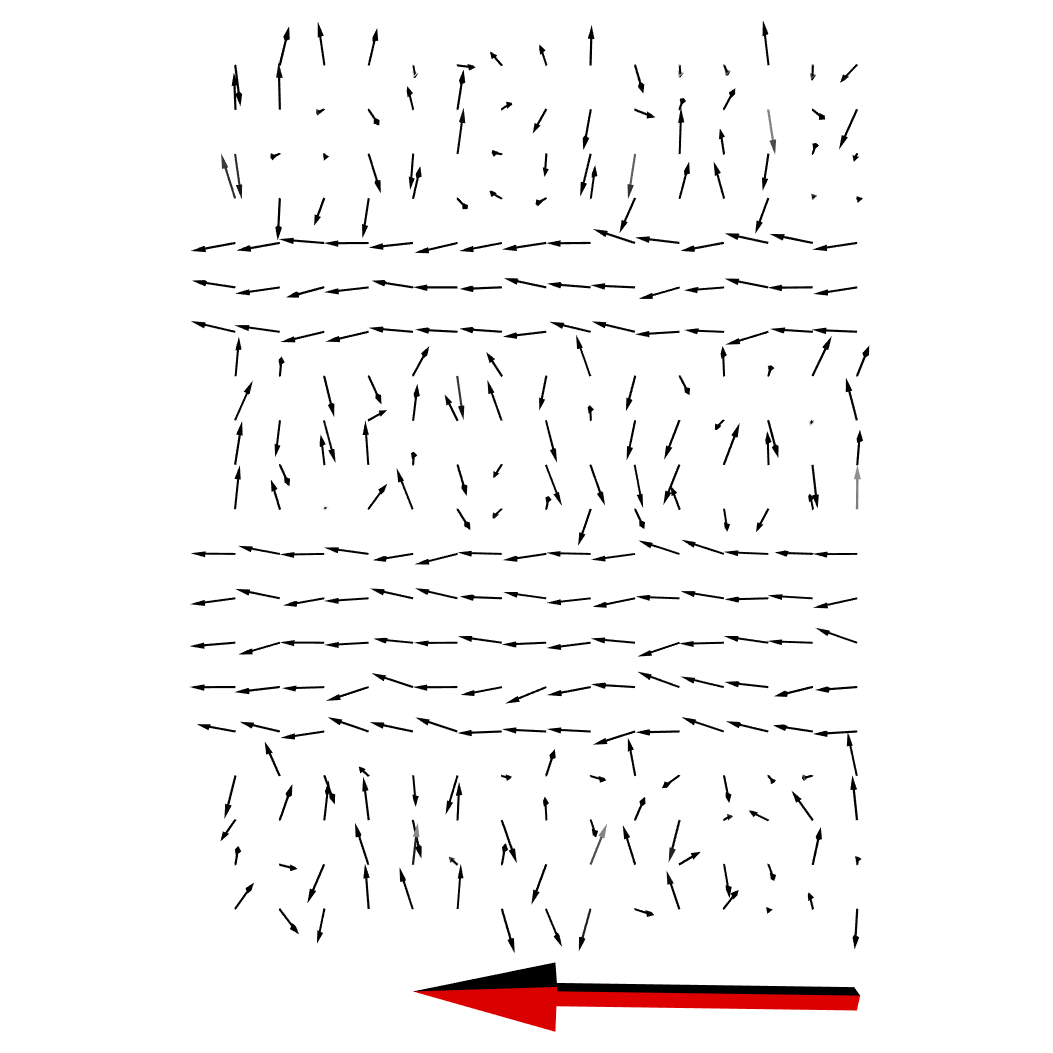}
    
   \caption{
    Original $\S$-valued image $\Iprime$.
    The direction vectors within the fibres are rotationally symmetric around $\mu$ (large red vector). 
    Their angular deviation from $\mu$ is uniformly drawn from $(0,\pi/8)$. 
    Background vectors are also rotationally symmetric around $\mu$. Their angular deviation from $\mu$ is uniformly drawn from $[3/8\pi,5/8\pi]$. 
    }
    \label{fig:orig_dilation_erosion}
\end{figure}

\begin{figure}
    \centering
    
    \begin{subfigure}[b]{0.29\linewidth}
        \includegraphics[width=\linewidth,angle =-90]{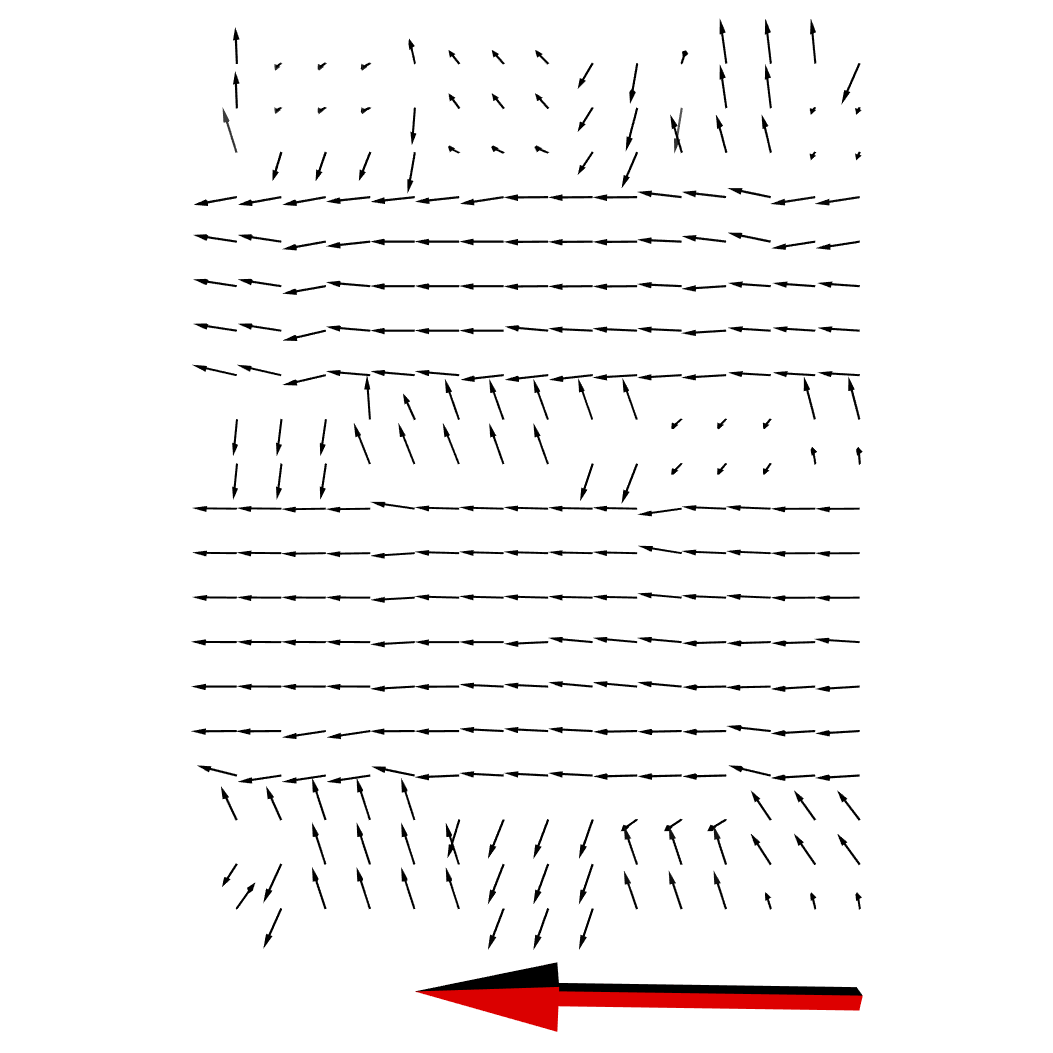}
        \caption{$\delta_B(\Iprime)$, $B$ is $3\times 3$}
    \end{subfigure}
    \begin{subfigure}[b]{0.29\linewidth}
        \includegraphics[width=\linewidth,angle =-90]{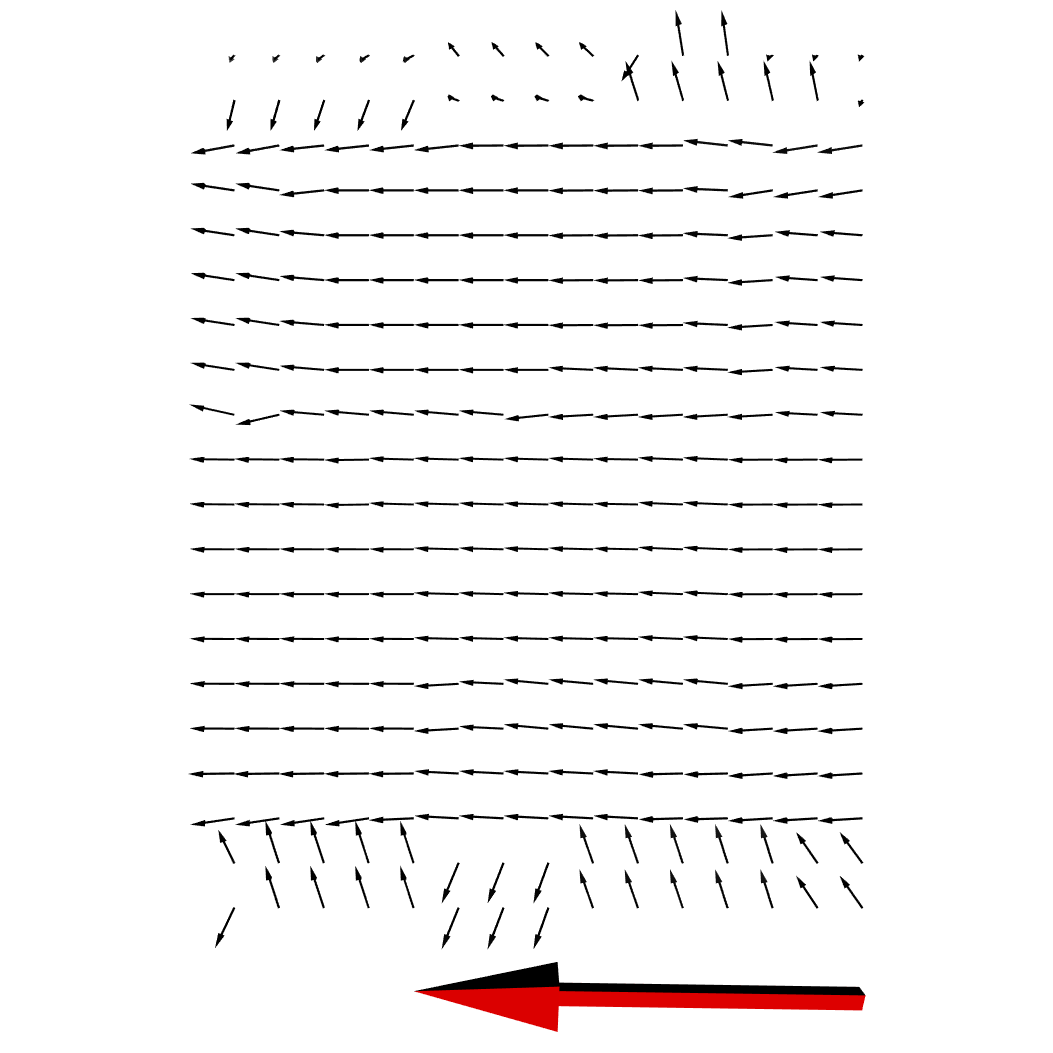}
        \caption{$\delta_B(\Iprime)$, $B$ is $5\times 5$}
    \end{subfigure}
    
    \begin{subfigure}[b]{0.29\linewidth}
        \includegraphics[width=\linewidth,angle =-90]{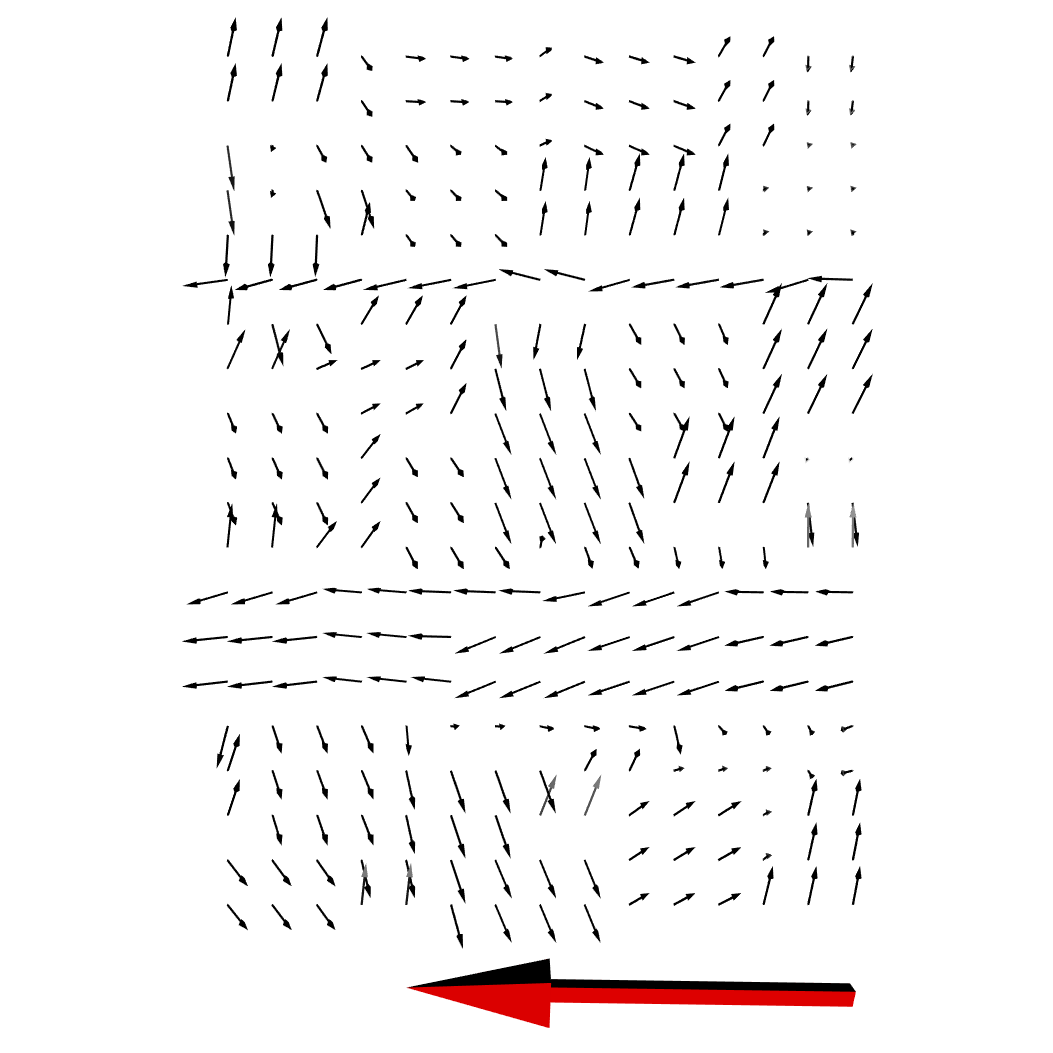}
        \caption{$\varepsilon_B(\Iprime)$, $B$ is $3\times 3$}
    \end{subfigure}
    \begin{subfigure}[b]{0.29\linewidth}
        \includegraphics[width=\linewidth,angle =-90]{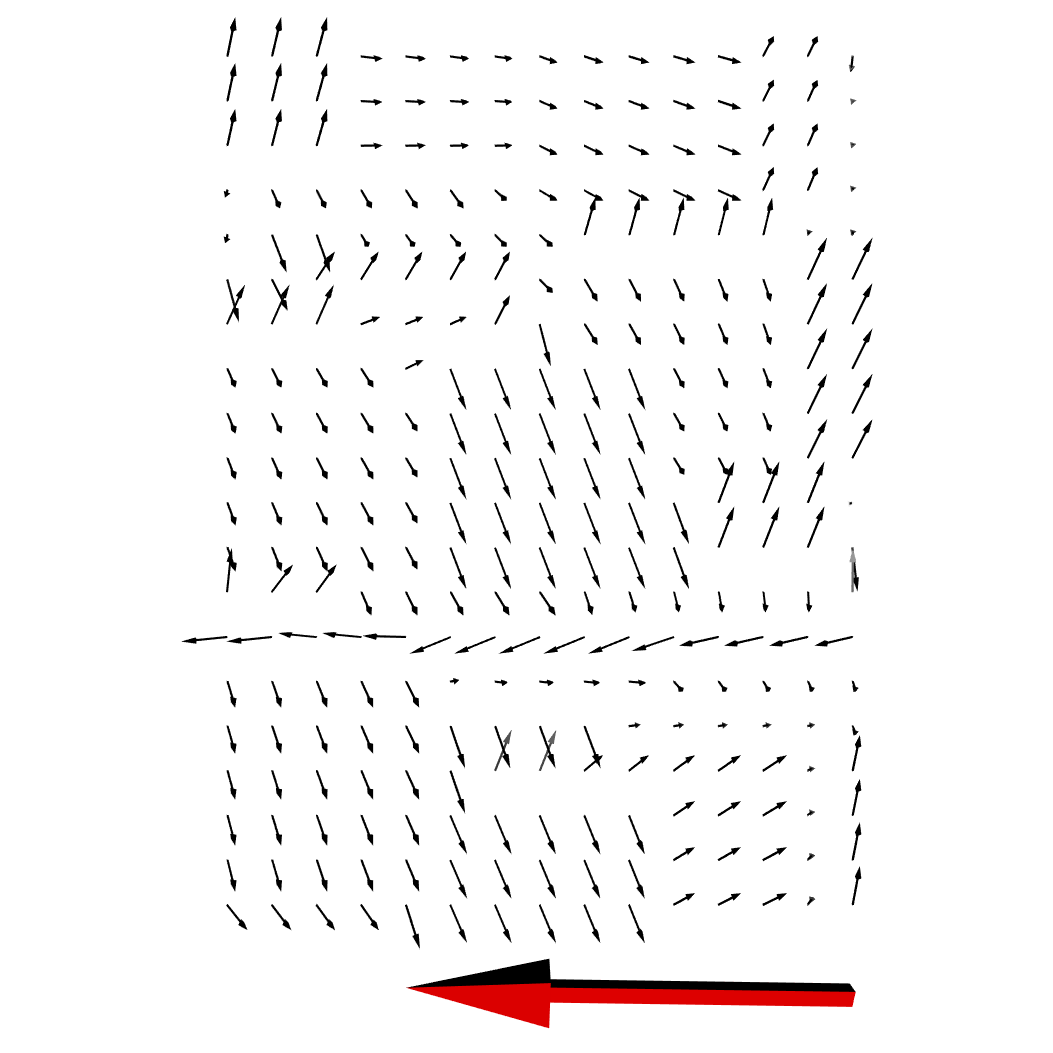}
        \caption{$\varepsilon_B(\Iprime)$, $B$ is $5\times 5$}
    \end{subfigure}

    \begin{subfigure}[b]{0.29\linewidth}
        \includegraphics[width=\linewidth,angle =-90]{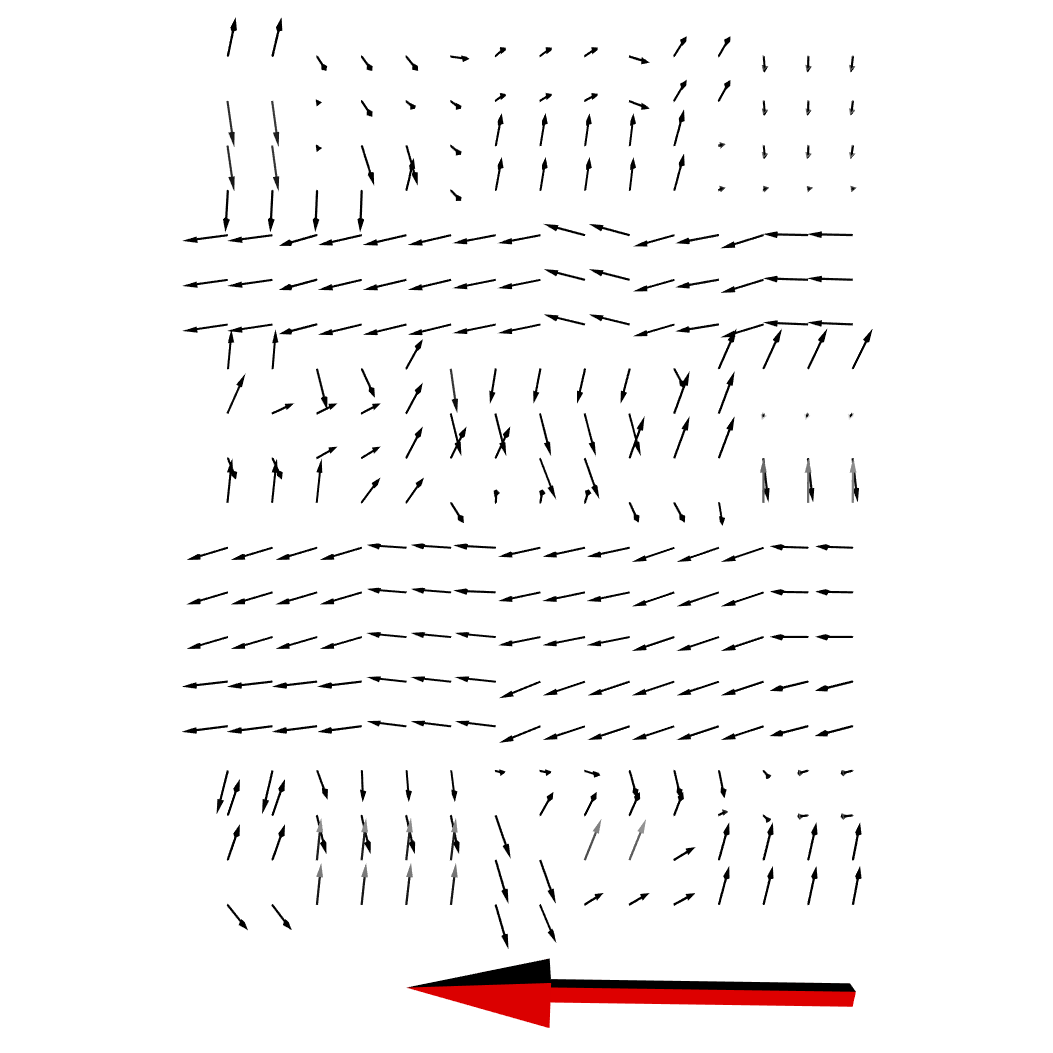}
        \caption{$\gamma_B(\Iprime)$, $B$ is $3\times 3$}
    \end{subfigure}
    \begin{subfigure}[b]{0.29\linewidth}
        \includegraphics[width=\linewidth,angle =-90]{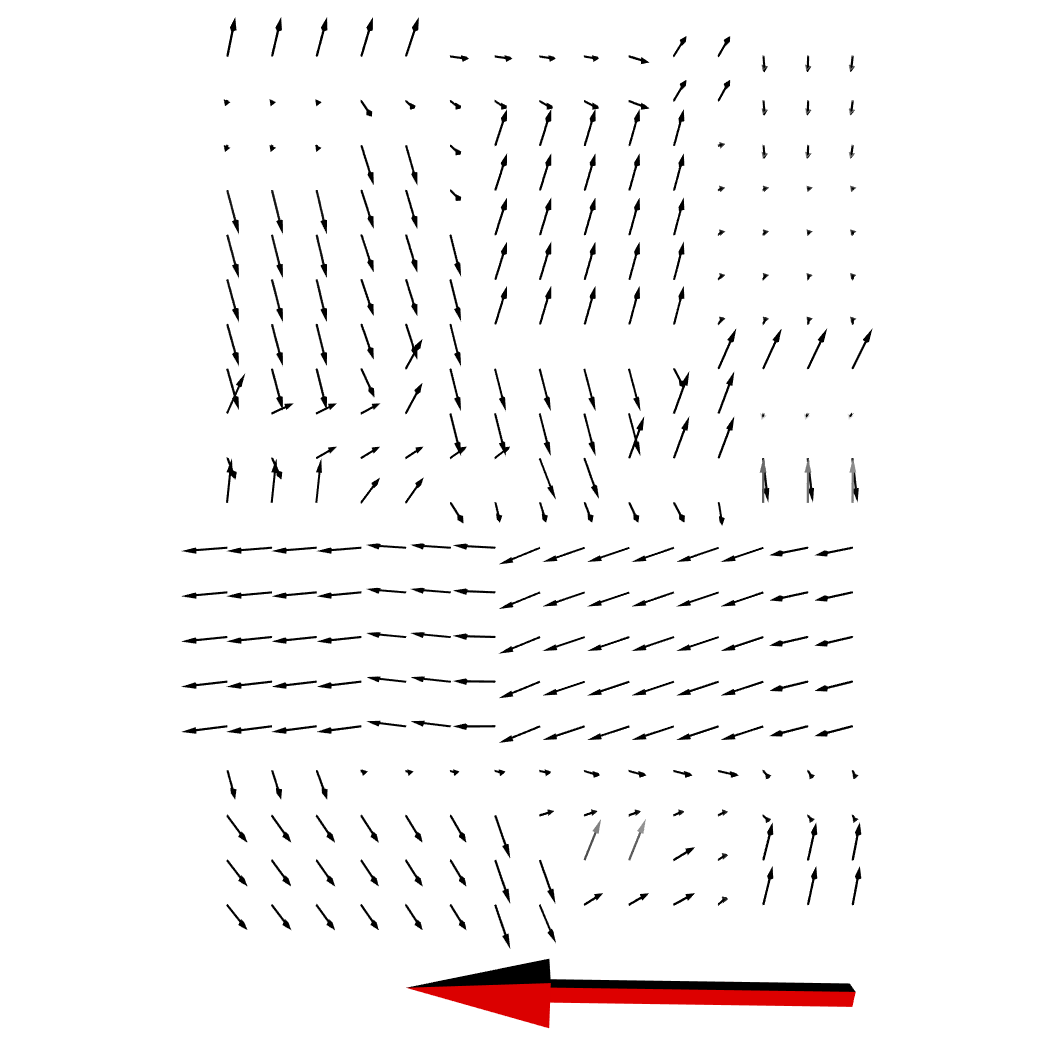}
        \caption{$\gamma_B(\Iprime)$, $B$ is $5\times 5$}
    \end{subfigure}

    \begin{subfigure}[b]{0.29\linewidth}
        \includegraphics[width=\linewidth,angle =-90]{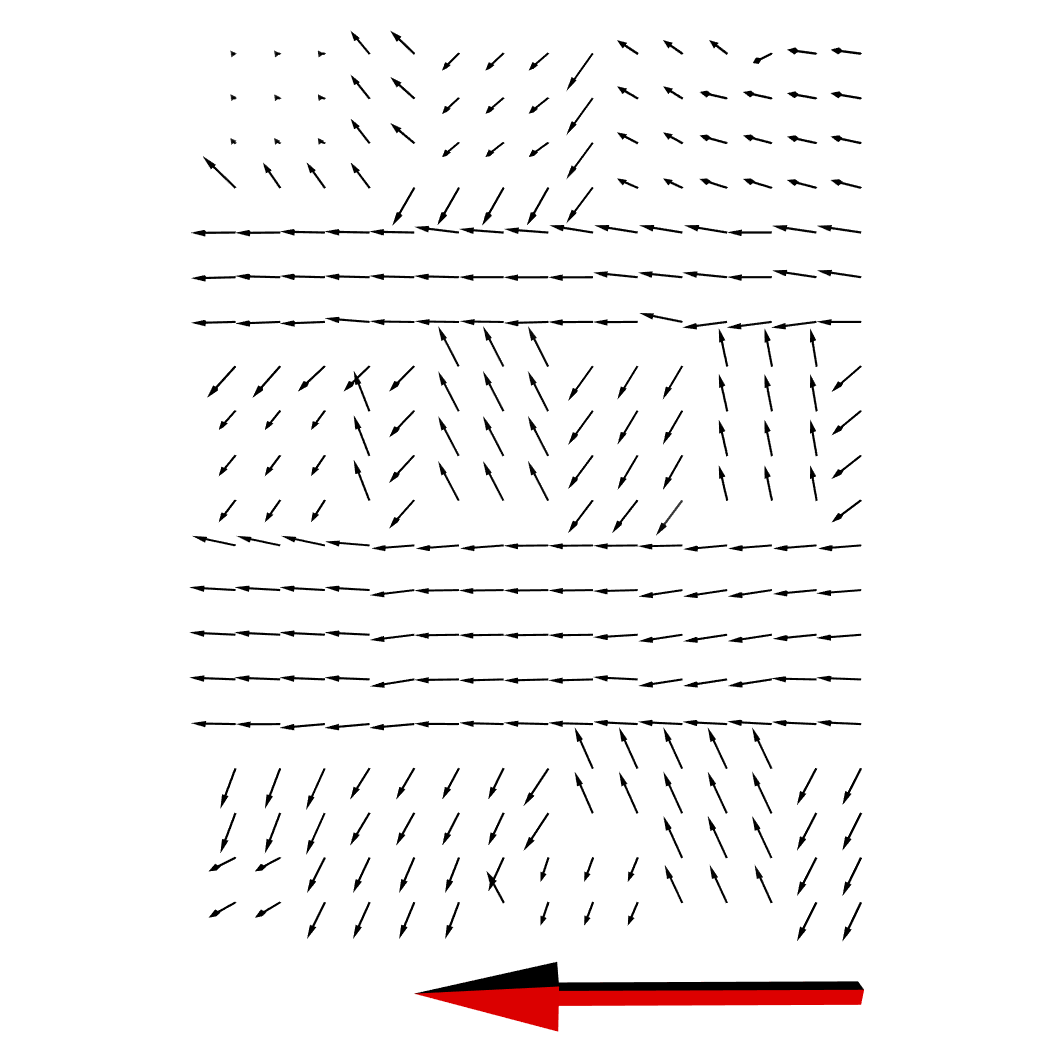}
        \caption{$\varphi_B(\Iprime)$, $B$ is $3\times 3$}
    \end{subfigure}
    \begin{subfigure}[b]{0.29\linewidth}
        \includegraphics[width=\linewidth,angle =-90]{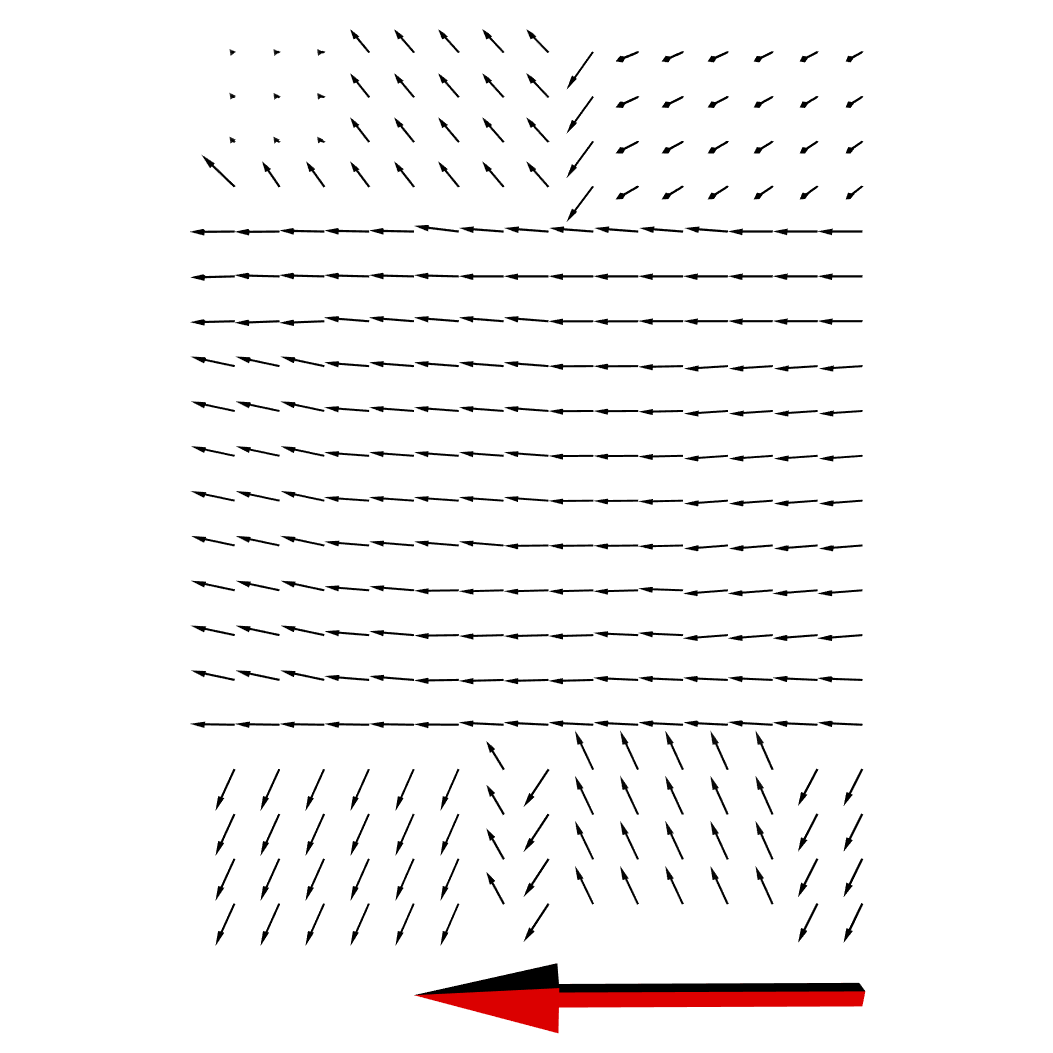}
        \caption{$\varphi_B(\Iprime)$, $B$ is $5\times 5$}
    \end{subfigure}
    \caption{
    Flat dilation $\delta_B$ (a,b),
    flat erosion $\varepsilon_B$ (c,d), 
    flat opening $\gamma_B$ (e,f)
    and flat closing $\varphi_B$ (g,h) 
    with $B$ a square as a structuring element. 
    The large red vector corresponds to $\mu$. 
    }
    \label{fig:flat_erosion_dilation_opening_closing}
\end{figure}

To investigate the newly defined flat morphological operators, we generate a 2-dimensional $\S$-valued image $\Iprime$ mimicking direction vectors obtained from two fibres on a homogeneous background.
That is, we assume that the image contains two objects consisting of vectors with a small angular deviation from $\mu$, see Figure \ref{fig:orig_dilation_erosion}. In the following, we will interpret those as foreground.
The background is formed by vectors directed along a plane perpendicular to $\mu$.
Flat dilation and flat erosion (Figure \ref{fig:flat_erosion_dilation_opening_closing}) show similar behaviour as their standard grey-scale counterparts.
Dilation expands objects directed along $\mu$. 
That is, the directions of background pixels that are close to the edge of these objects are rotated towards $\mu$. 
An erosion shrinks objects, i.e., 
object vectors at the edge are rotated so that they are assigned to the image background.

In a similar manner, the flat opening removes foreground objects that are smaller than the structuring element. 
Figure \ref{fig:flat_erosion_dilation_opening_closing} reveals that object vectors in an object of size smaller than the structuring element $B$ are rotated such that we assign them to the image background.
The flat closing removes small holes in the foreground. 
Vectors with large angular deviation from $\mu$ within a background object of size smaller than $B$ are rotated to become part of the image foreground.

\begin{figure}
    \centering
    \begin{subfigure}[b]{0.48\linewidth}
        \includegraphics[width=\linewidth,angle =-90]{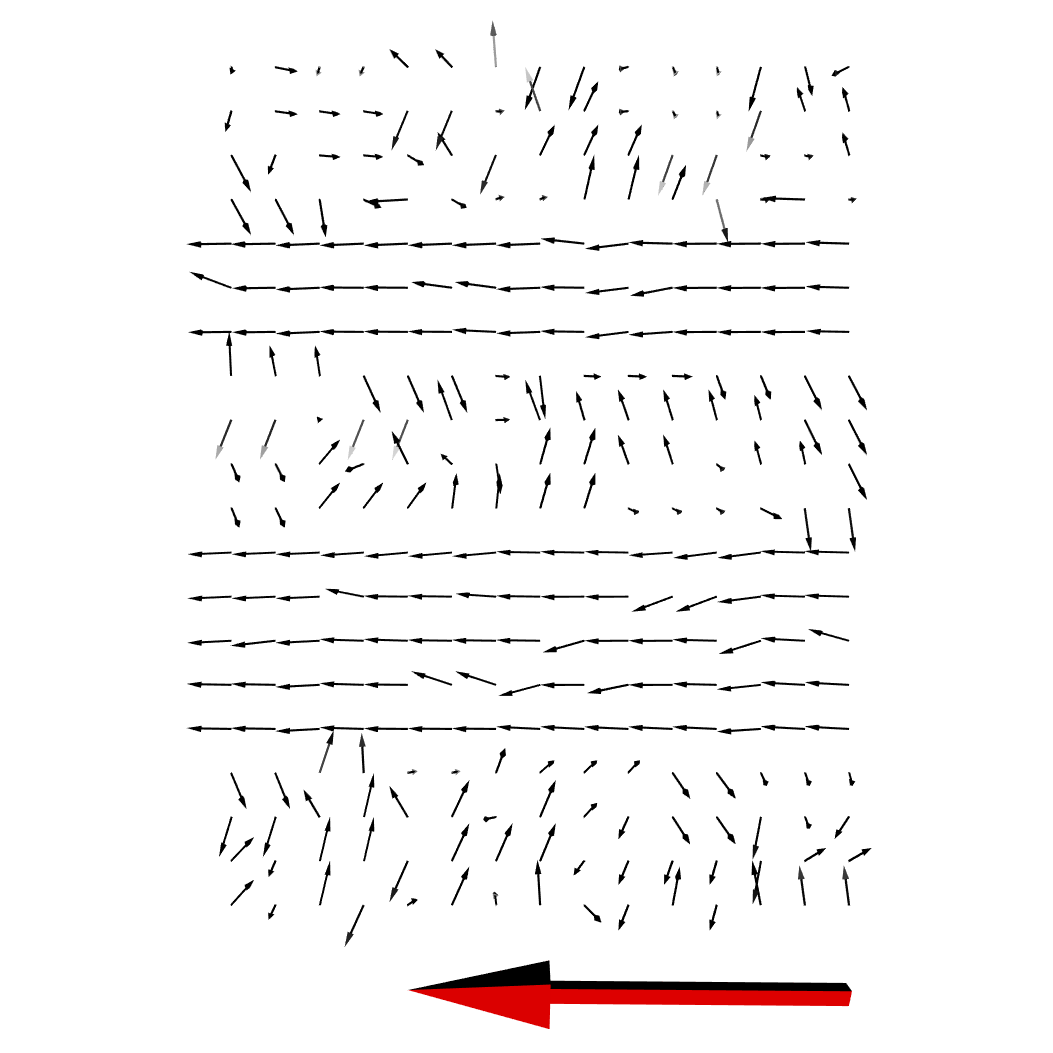}
        \caption{$\shockfilter_B(\Iprime)$, $B$ is $3\times 3$}
    \end{subfigure}
    \begin{subfigure}[b]{0.48\linewidth}
        \includegraphics[width=\linewidth,angle =-90]{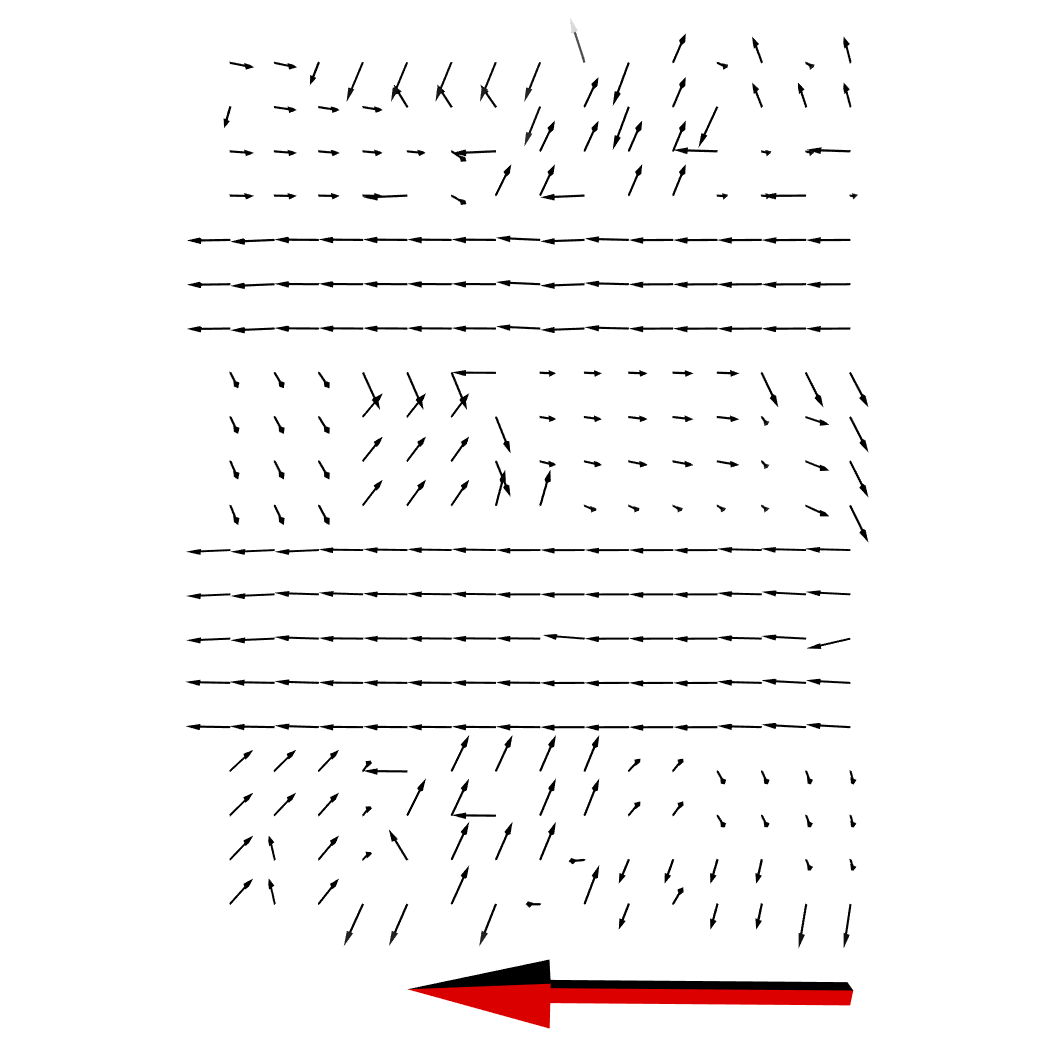}
        \caption{$\shockfilter_B(\Iprime)$, $B$ is $5\times 5$}
    \end{subfigure}

    \caption{
    Flat shock filter $\shockfilter_B(\Iprime)$ with $B$ a square as structuring element. 
    }
    \label{fig:flat_shock_filter}
\end{figure}

The results of shock filtering of $\Iprime$ are shown in Figure \ref{fig:flat_shock_filter}.
The edges between the two objects and the background pixels are enhanced.
\subsection{Pseudo morphological multi-scale operators}
\begin{figure}
	\centering
	\begin{subfigure}[b]{0.39\linewidth}
		\includegraphics[width=\linewidth,angle =-90]{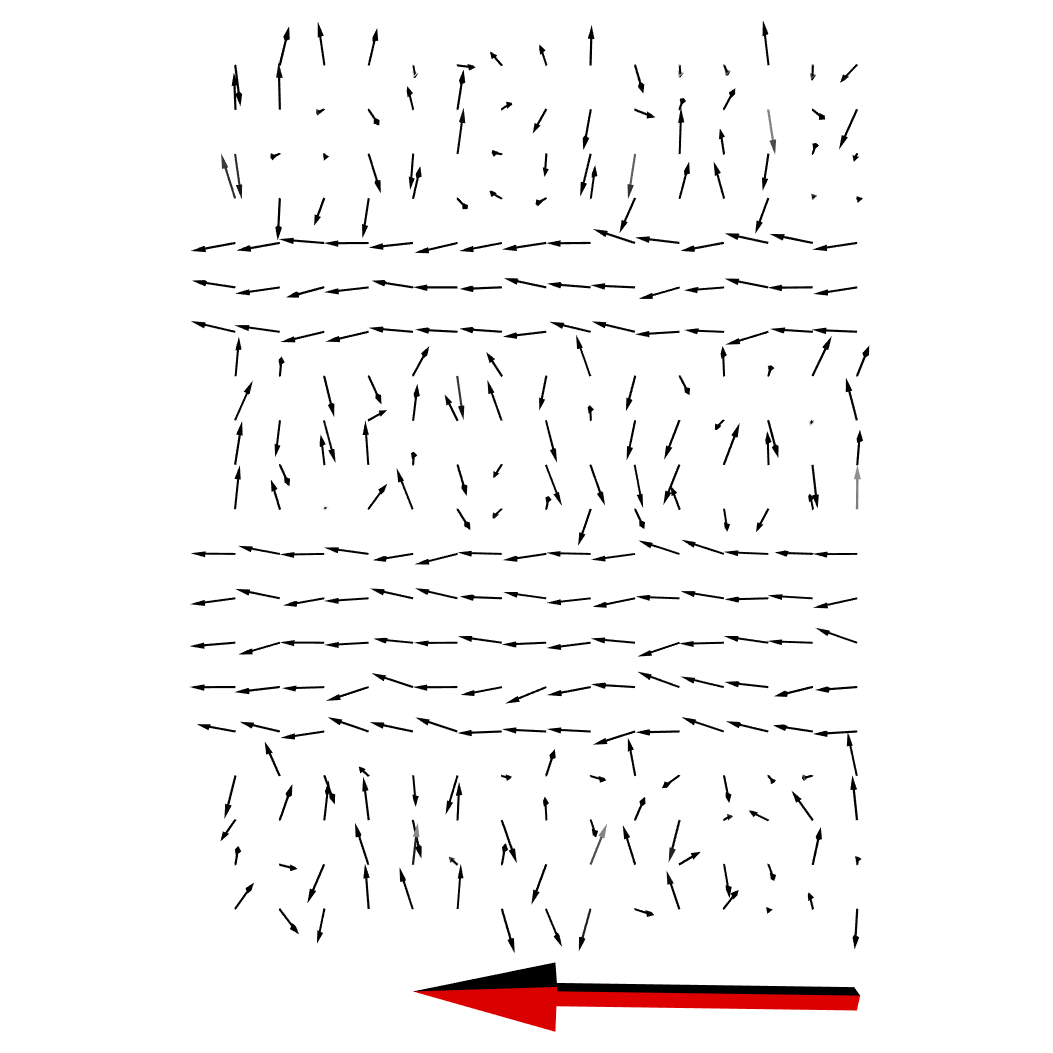}
		\caption{$t=0.1$ (equals input image)}
	\end{subfigure}
	\begin{subfigure}[b]{0.39\linewidth}
		\includegraphics[width=\linewidth,angle =-90]{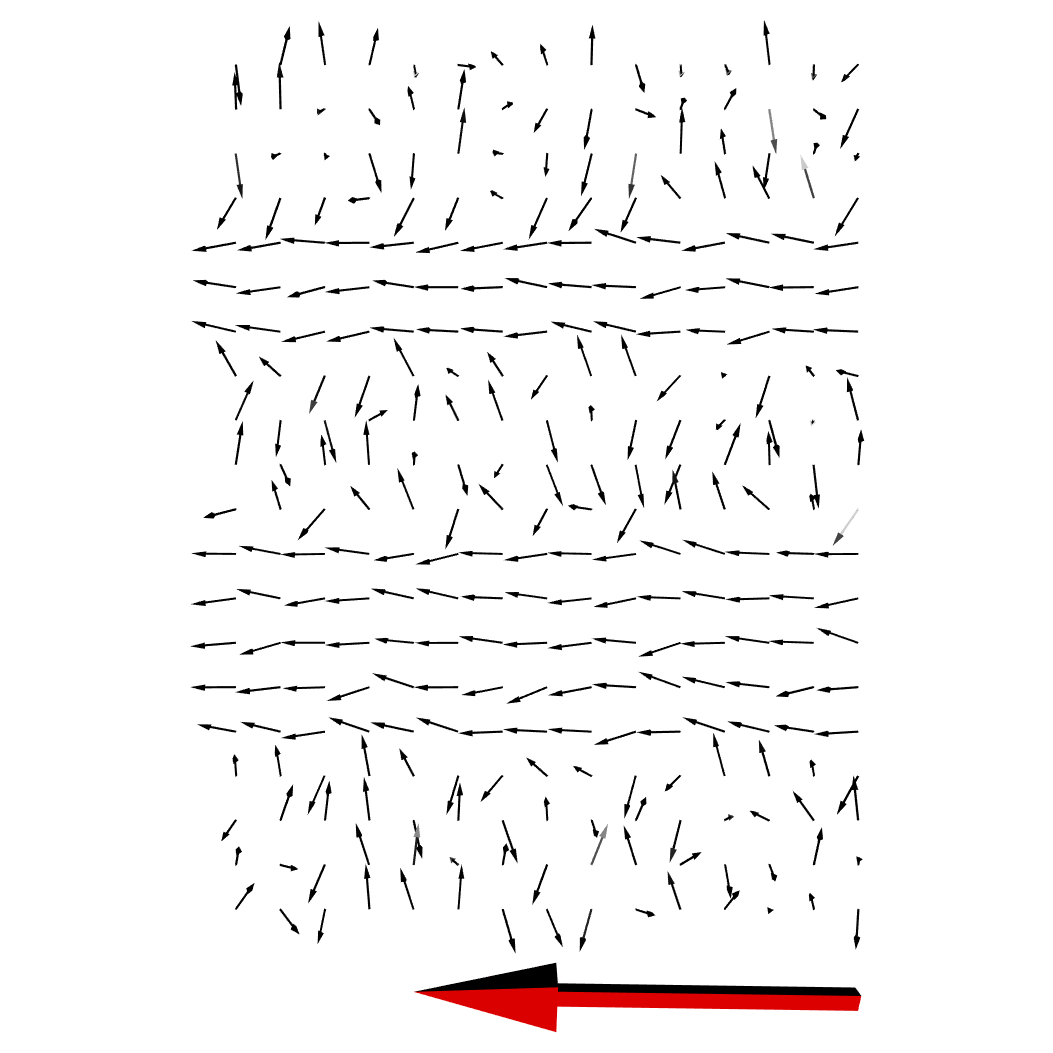}
		\caption{$t=0.5$}
	\end{subfigure}
	
	\begin{subfigure}[b]{0.39\linewidth}
		\includegraphics[width=\linewidth,angle =-90]{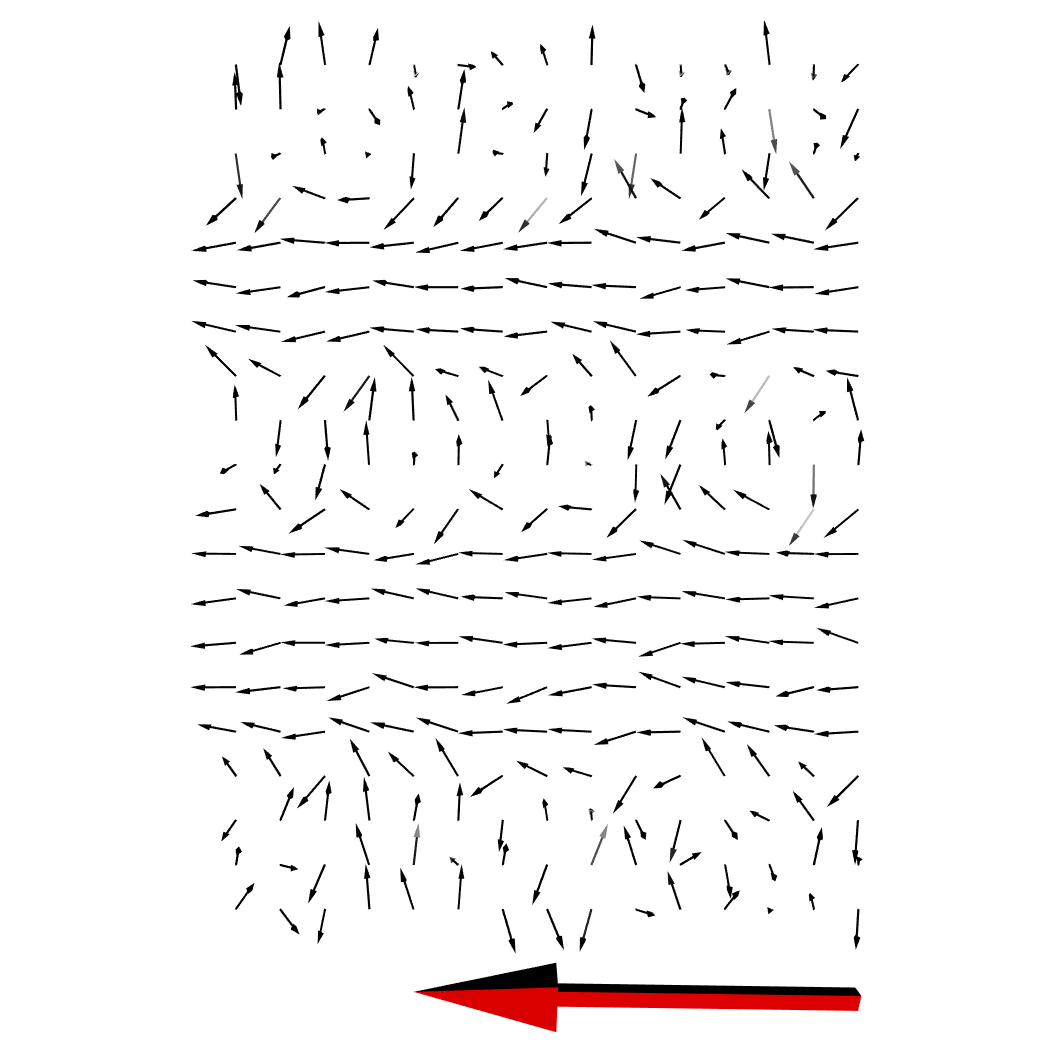}
		\caption{$t=0.7$}
	\end{subfigure}
	\begin{subfigure}[b]{0.39\linewidth}
		\includegraphics[width=\linewidth,angle =-90]{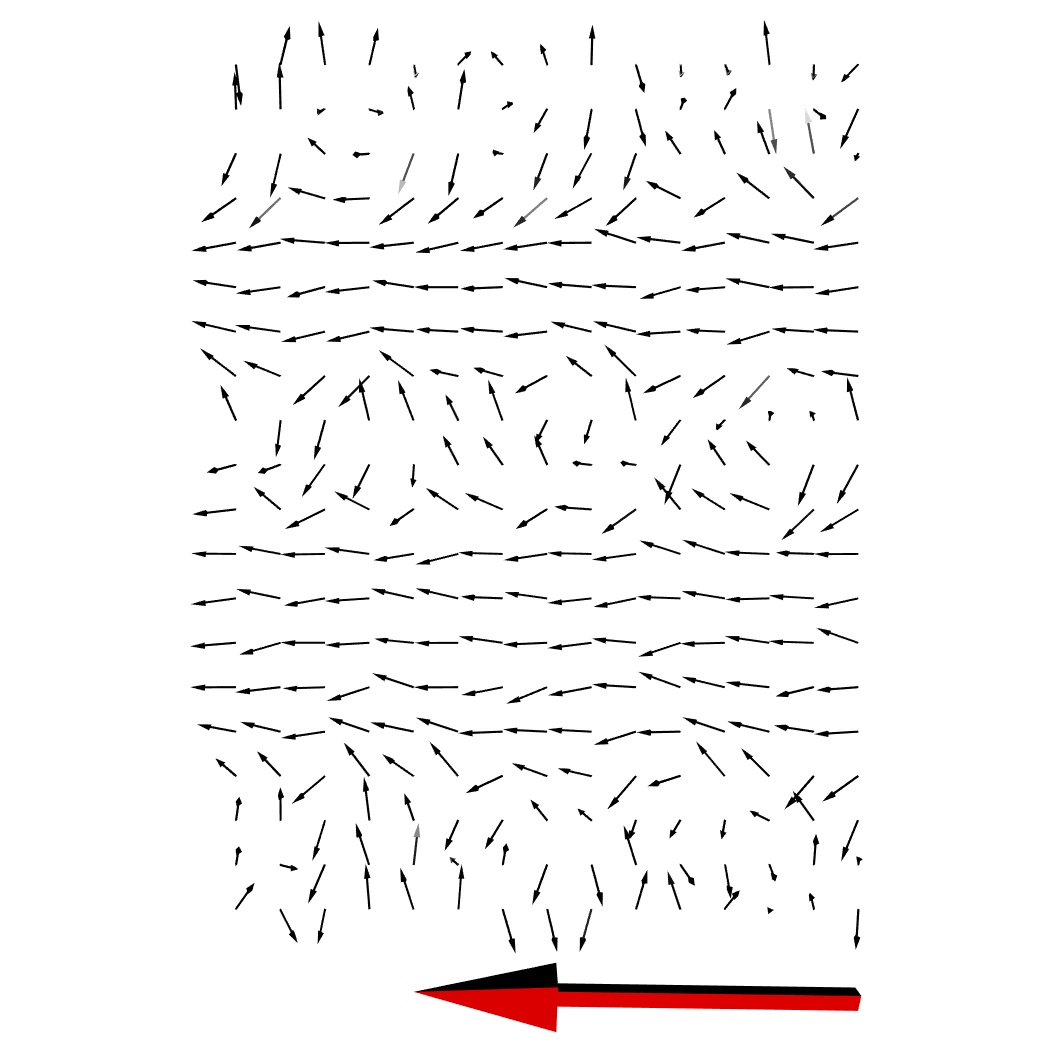}
		\caption{$t=0.9$}
	\end{subfigure}
	
	\begin{subfigure}[b]{0.39\linewidth}
		\includegraphics[width=\linewidth,angle =-90]{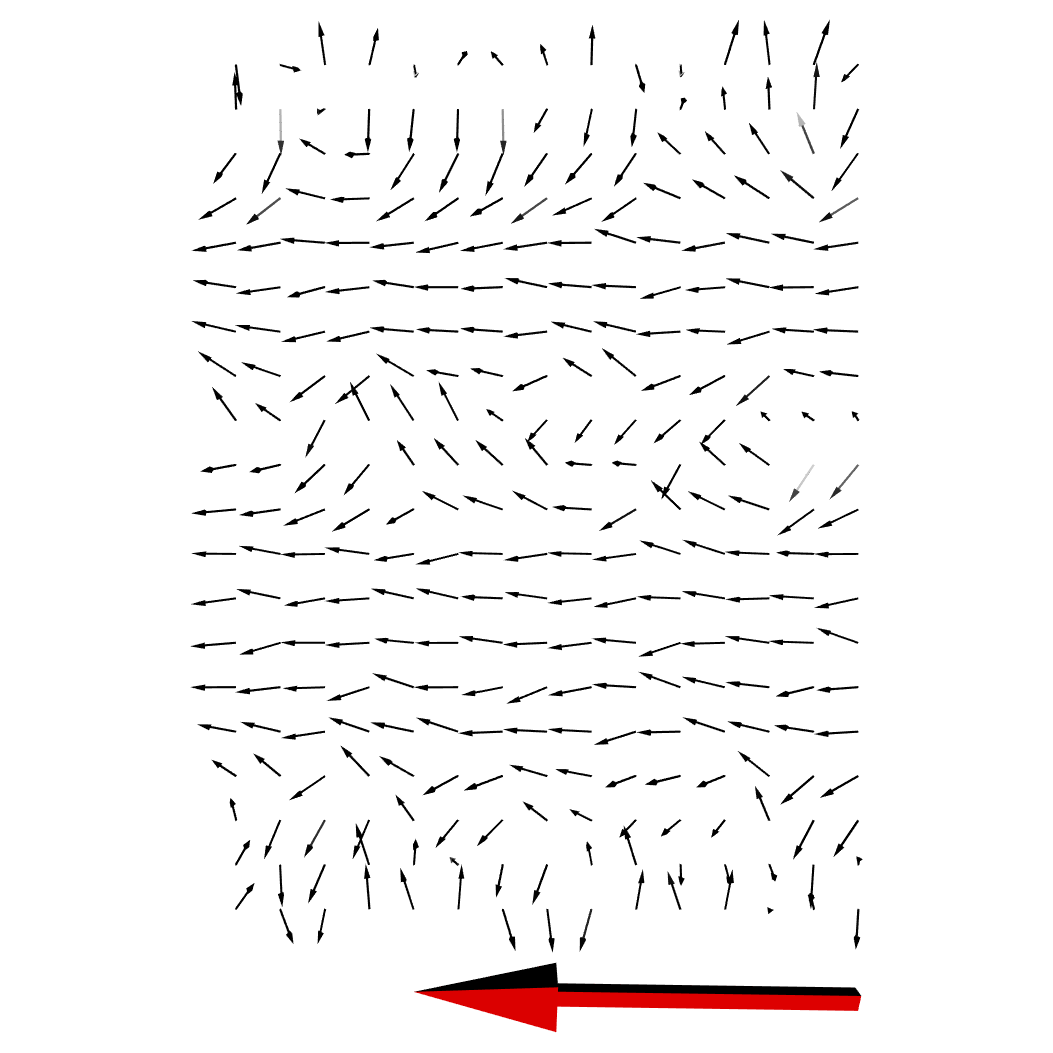}
		\caption{$t=1.1$}
	\end{subfigure}
	\begin{subfigure}[b]{0.39\linewidth}
		\includegraphics[width=\linewidth,angle =-90]{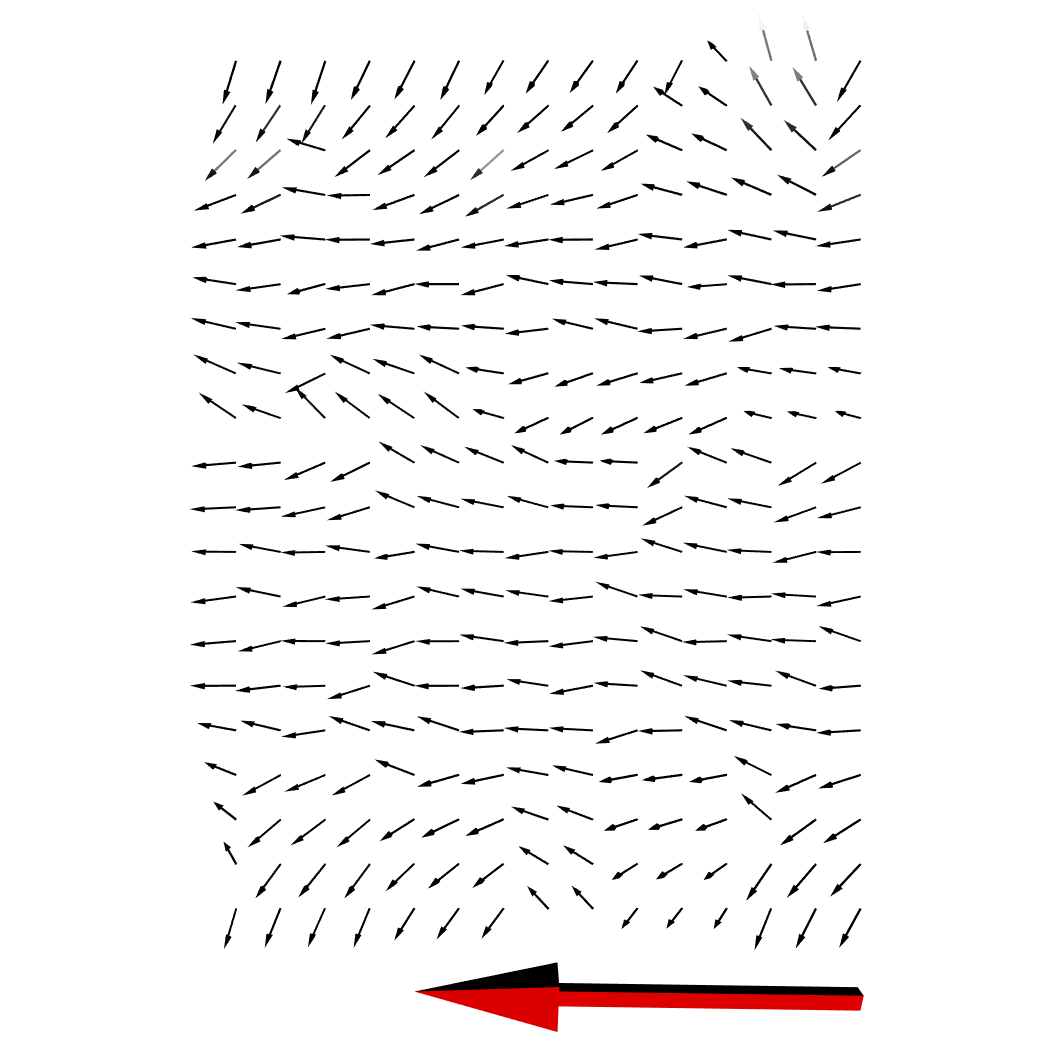}
		\caption{$t=2$}
	\end{subfigure}
	
	\caption{$\delta_\btstwo(\Iprime)$ at scales $t=0.1,0.5,0.7,0.9,1.1,2$. The large red vector corresponds to $\mu$.}
	\label{fig:scale_space_dilation}
\end{figure}

\begin{figure}
	\centering
	\begin{subfigure}[b]{0.39\linewidth}
		\includegraphics[width=\linewidth,angle =-90]{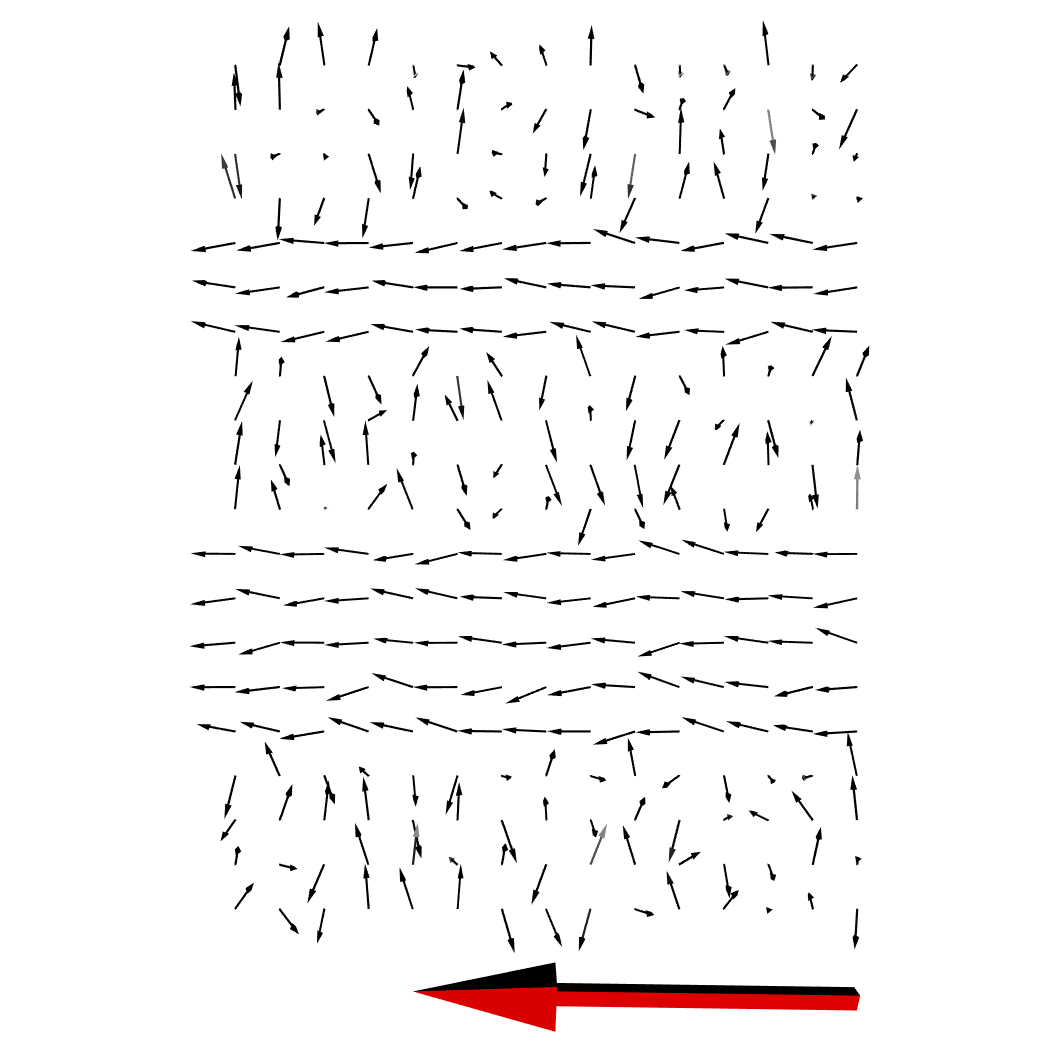}
		\caption{$t=0.1$}
	\end{subfigure}
	\begin{subfigure}[b]{0.39\linewidth}
		\includegraphics[width=\linewidth,angle =-90]{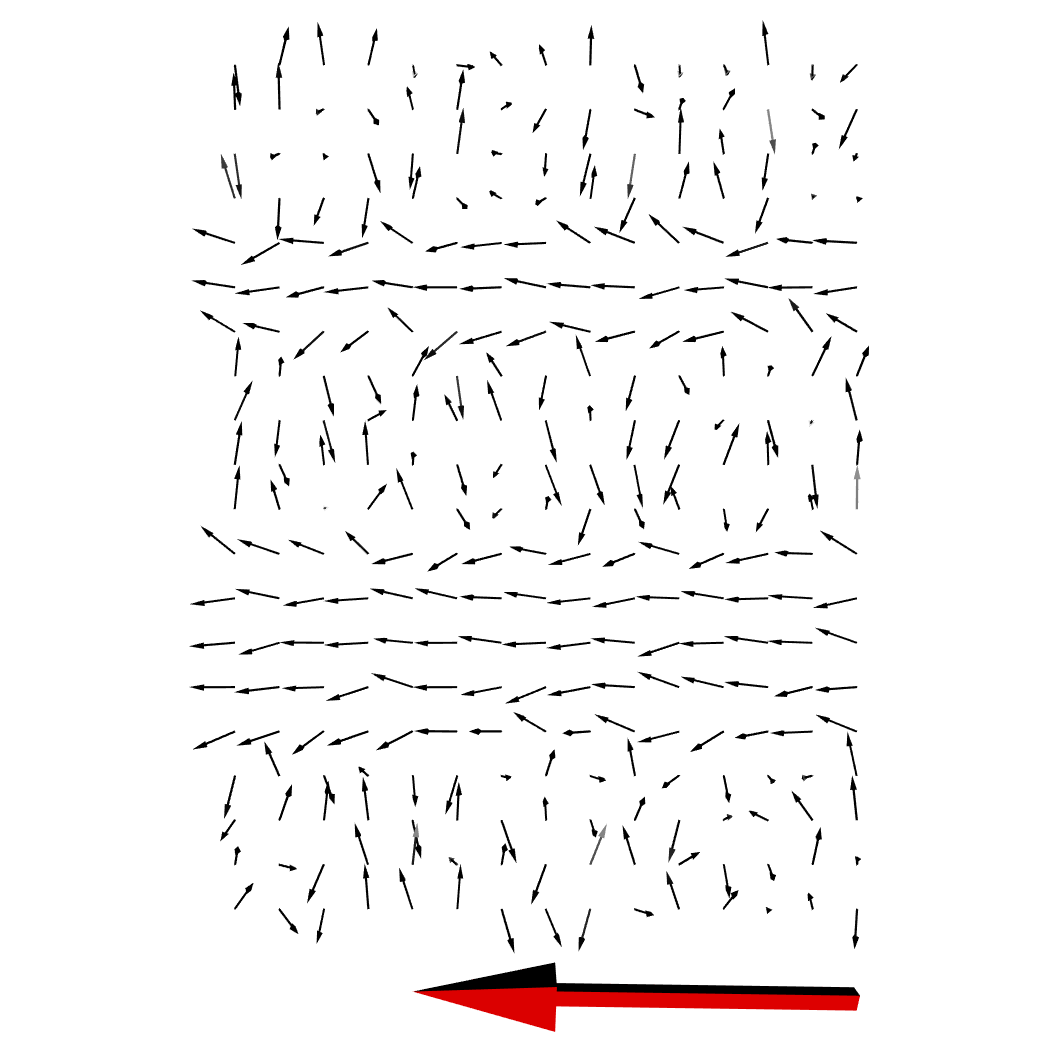}
		\caption{$t=0.5$}
	\end{subfigure}
	
	\begin{subfigure}[b]{0.39\linewidth}
		\includegraphics[width=\linewidth,angle =-90]{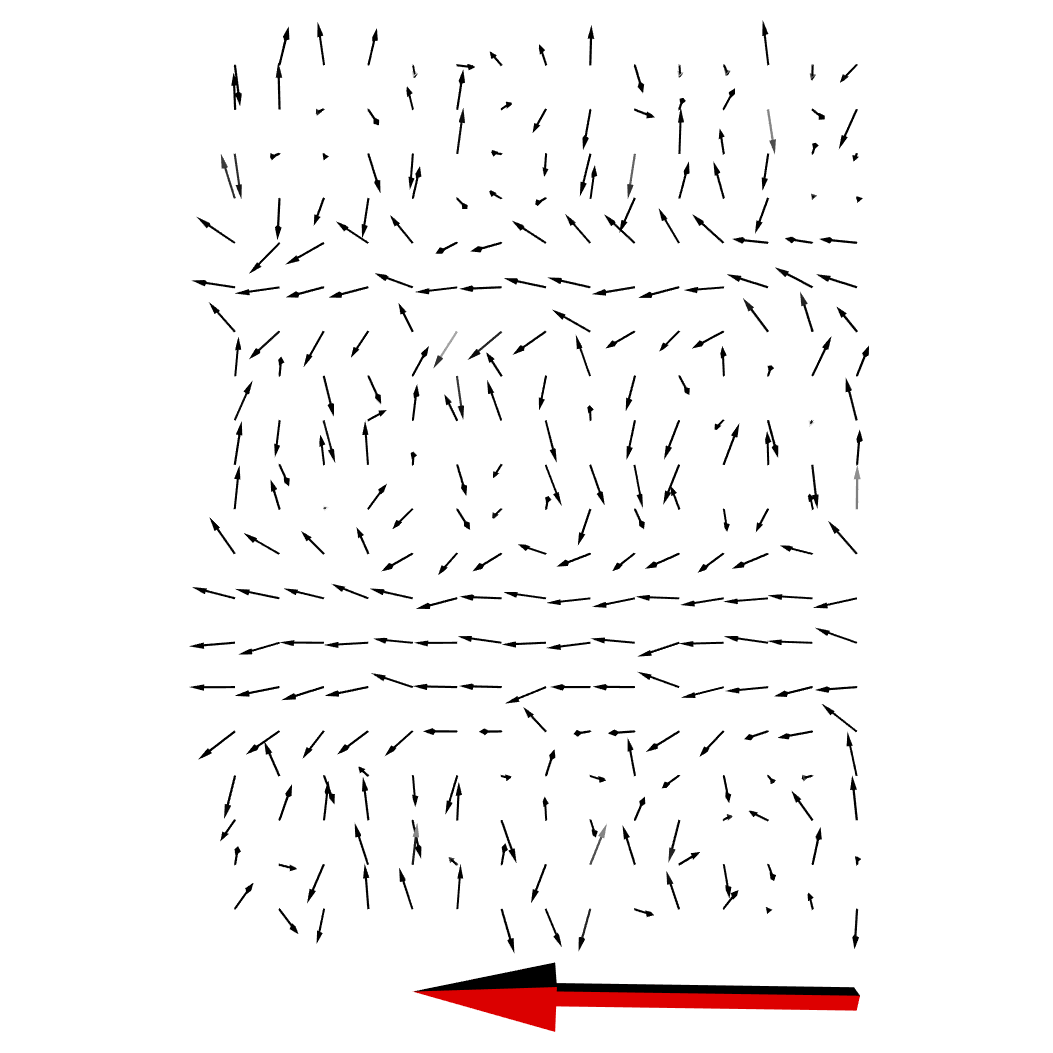}
		\caption{$t=0.7$}
	\end{subfigure}
	\begin{subfigure}[b]{0.39\linewidth}
		\includegraphics[width=\linewidth,angle =-90]{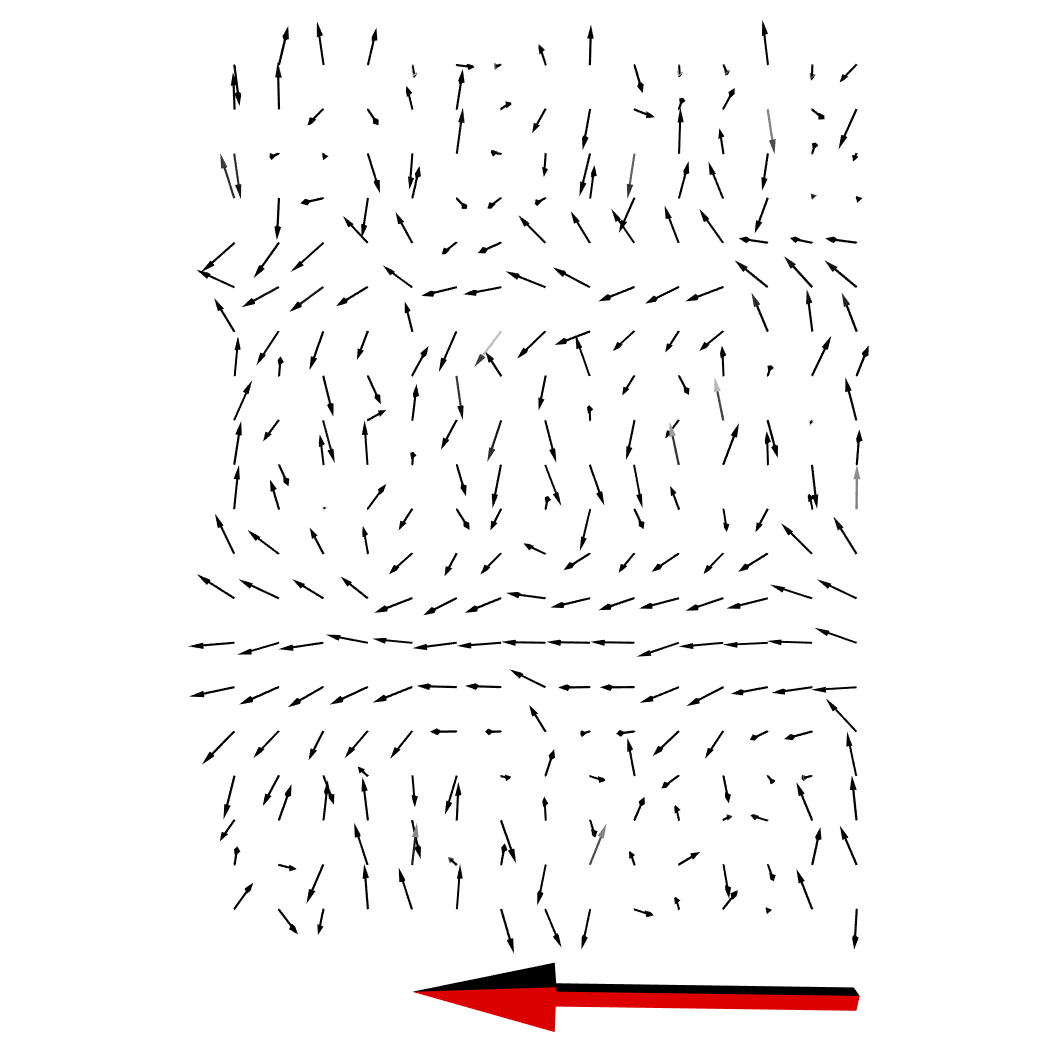}
		\caption{$t=0.9$}
	\end{subfigure}
	
	\begin{subfigure}[b]{0.39\linewidth}
		\includegraphics[width=\linewidth,angle =-90]{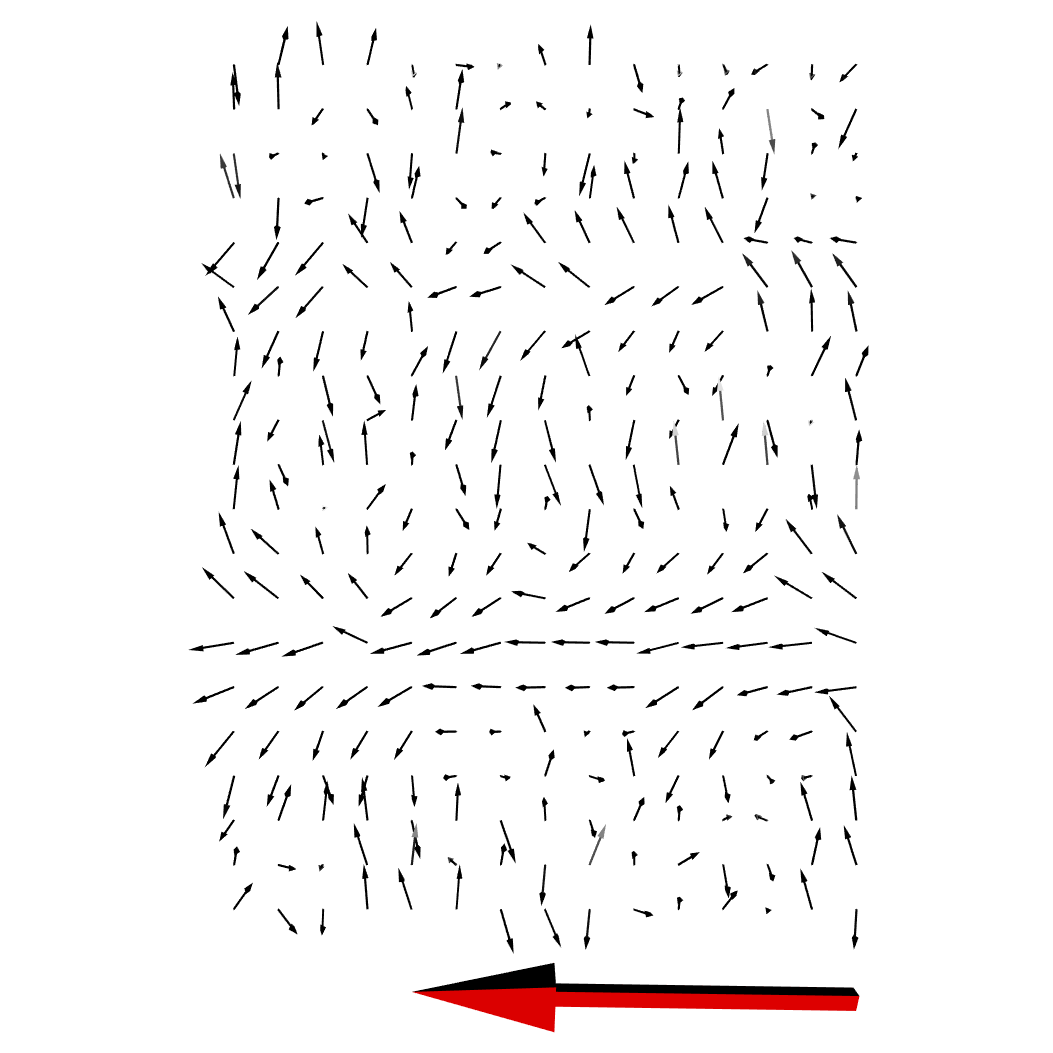}
		\caption{$t=1.1$}
	\end{subfigure}
	\begin{subfigure}[b]{0.39\linewidth}
		\includegraphics[width=\linewidth,angle =-90]{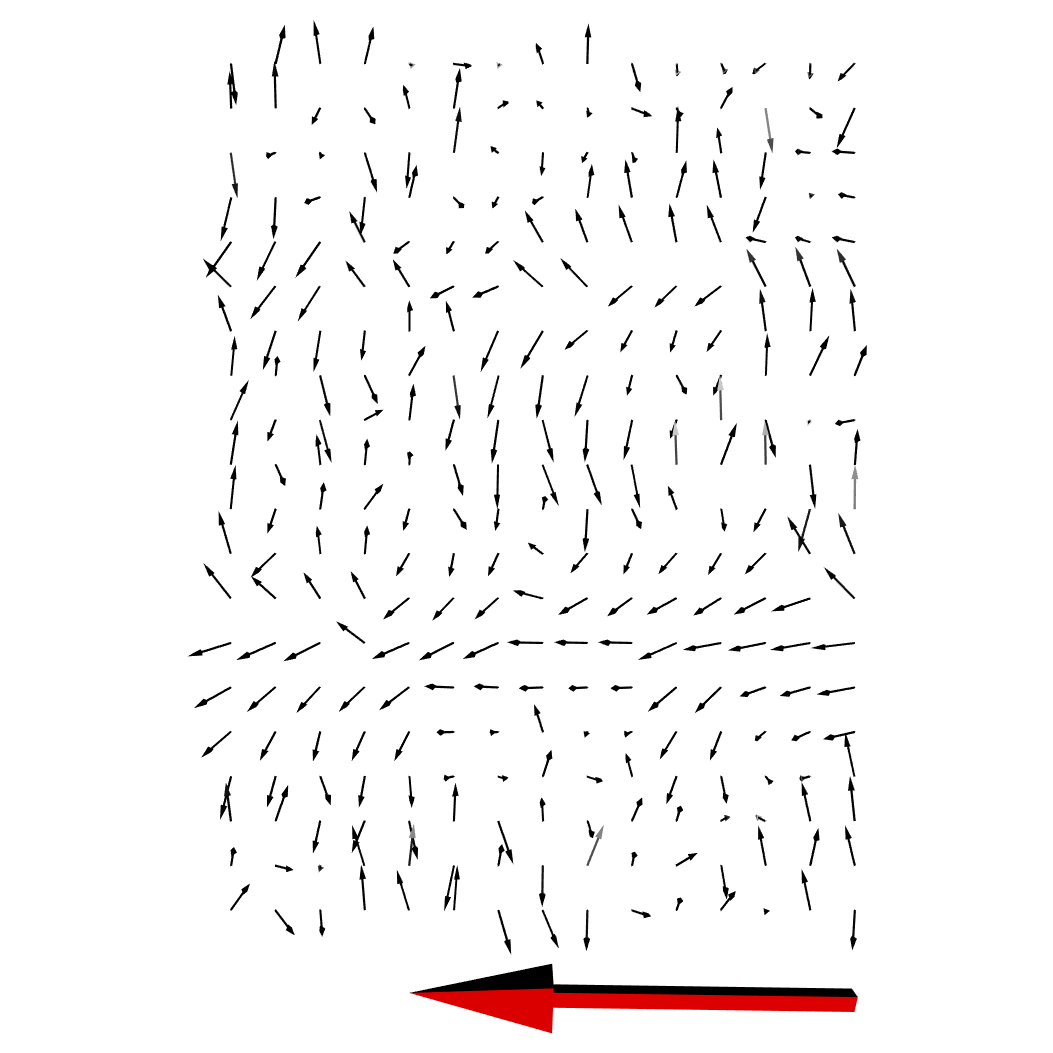}
		\caption{$t=2$}
	\end{subfigure}
	
	\caption{$\varepsilon_\btstwo(\Iprime)$ at scales $t=0.1,0.5,0.7,0.9,1.1,2$. The large red vector corresponds to $\mu$.}
	\label{fig:scale_space_erosion}
\end{figure}

We illustrate the multi-scale dilation and erosion on a directional image $\Iprime$. 
$\delta_\btstwo(\Iprime)$ and $\varepsilon_\btstwo(\Iprime)$ behave like their grey-scale counterparts:
Figure \ref{fig:scale_space_dilation} shows that $\delta_\btstwo(\Iprime)$ enlarges objects in the foreground and Figure \ref{fig:scale_space_erosion} shows that $\varepsilon_\btstwo(\Iprime)$ shrinks objects in the foreground. 
Of course, the enlargement and shrinkage depend on the continuous scale parameter $t$.





\subsection{Segmentation of aligned and misaligned fibre regions of glass fibre reinforced polymers}

Glass fibre reinforced polymers (GFRP) are frequently used in lightweight construction of, e.g.,  cars, planes, and wind power plants. In civil engineering, they have emerged as an alternative to steel reinforcement as they do not corrode in aggressive conditions such as marine environments \cite{Kappenthuler2021}.
GFRP materials are anisotropic and characterised by high tensile strength in the direction of the reinforcing fibres.
When manufactured by injection moulding, the fibres follow the injection direction.
However, the fibre orientation in a central layer differs depending on the production parameters.
Analysing fibre orientation and quantifying misalignment, for instance by directional filtering \cite{Tin22}, may help to optimise the production process and to improve GFRP.

\begin{figure}
     \centering
     \begin{subfigure}[b]{\linewidth}
     \centering
     \includegraphics[width=0.81\linewidth]{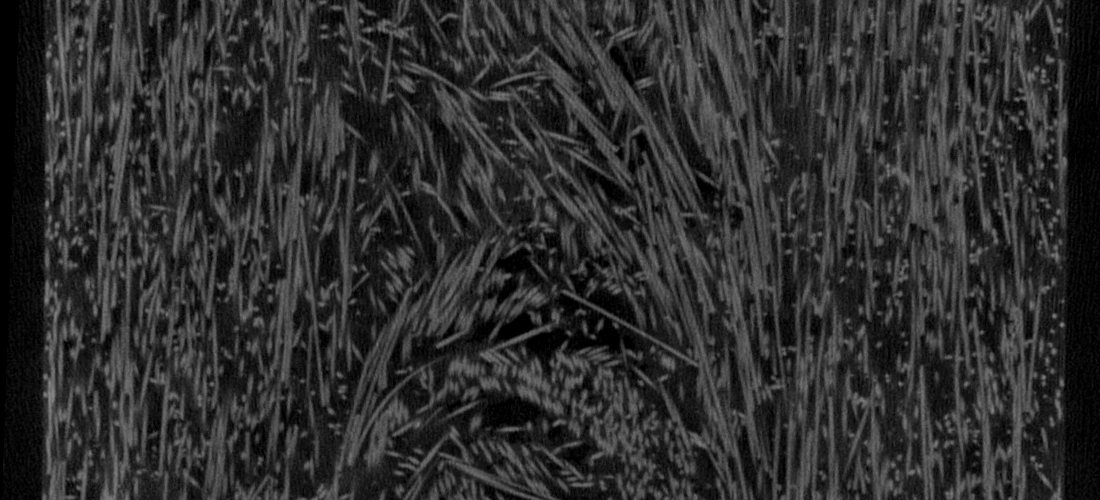}
     \caption{2D slice of the original reconstructed 3D $\mu$CT image. The image size is 1100$\times$500 pixels with a pixel spacing of 4 $\mu$m.}\label{fig:Original_CT_image_Tin}
     \end{subfigure}

     \begin{subfigure}[b]{\linewidth}
     \centering
     \includegraphics[width=0.4\linewidth]{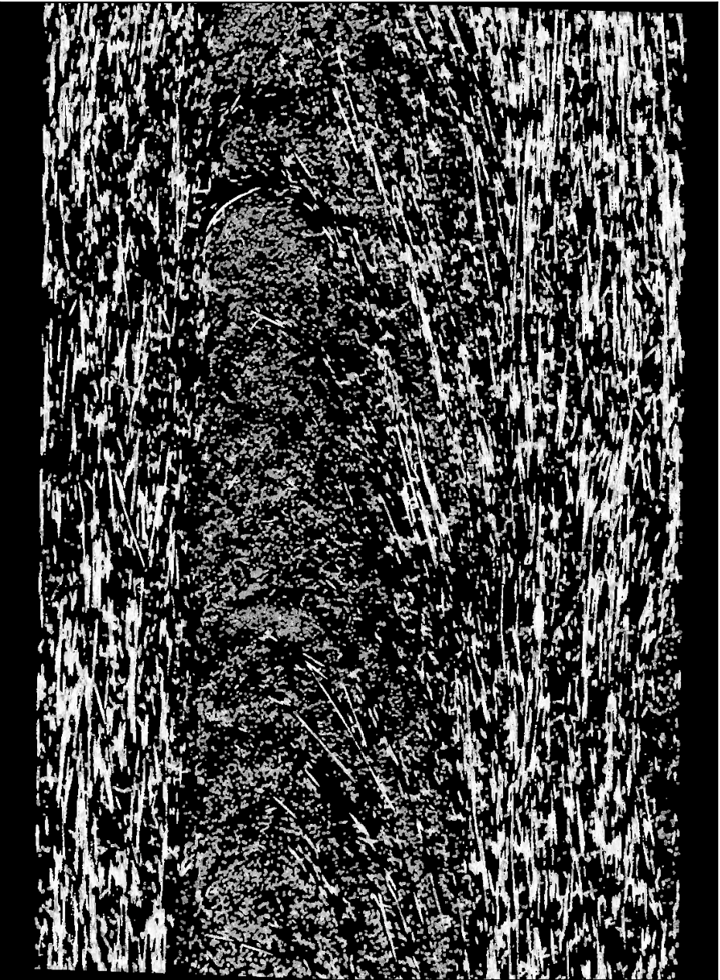}     \includegraphics[width=0.4\linewidth]{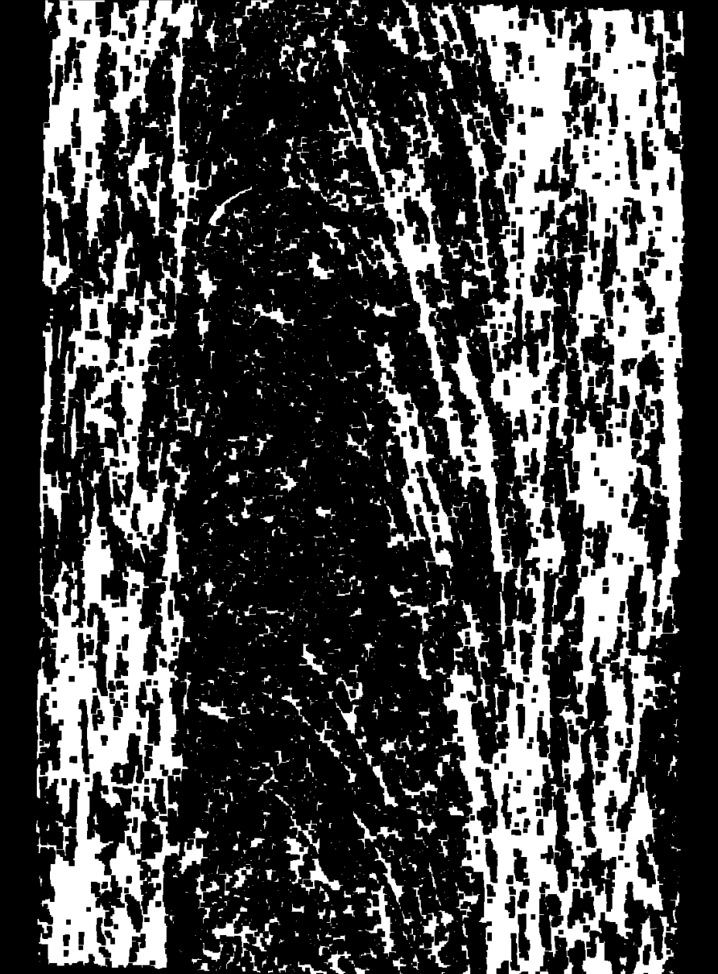}
     \caption{Left: 2D slice of a 3D grey scale image encoding the angular deviation of the fibre direction in each pixel from $\mu=(0,1,0)$.  The left and right side are aligned with $\mu$ (correctly aligned). In the inner region, the fibres are misaligned. 
     Right: Closing of the directional image (with $B$ a $7\times 7\times 7$ cube) and binarisation of the resulting angular deviation image (threshold 1.5) highlight the correctly aligned areas. The image size is 1100 $\times$ 1500 pixels.}\label{fig:closing_7x7x7_718_Bin_1.5}
     \end{subfigure}
    

    \caption{Sectional $\mu$CT image of GFRP produced and imaged by the Leibniz Institute for Composite Materials (IVW) in Kaiserslautern (Germany), the corresponding original and closed directional image.}
 \end{figure}

We examine a $\mu$CT image of GFRP produced and imaged by the Leibniz Institute for Composite Materials (IVW) in Kaiserslautern (Germany) and discussed in \cite{Wirjadi2014}. 
See Figure \ref{fig:Original_CT_image_Tin} for an illustration.
Pixel-wise fibre orientations can be determined by using partial second derivatives as described in \cite{ImageAnalStereol1489} resulting in an $\S$-valued image.


Our goal is to segment the region aligned along the injection direction, as desired, and the misaligned fibre region. 
Figure \ref{fig:closing_7x7x7_718_Bin_1.5} shows a grey value encoding of the scalar product between the pixelwise fibre directions and $\mu$. 
The correctly aligned regions are at the left and right borders while fibres in the central part are misaligned. The higher the value (brighter the image pixel), the more aligned the vectors are along $\mu$. 
Simply thresholding this image results in a very noisy segmentation. Therefore, the directional image is smoothed prior to the segmentation by a morphological closing along $\mu=(0,1,0)$ with $B$ a $7\times 7\times 7$ cube. 

We choose a threshold value of 1.5 to segment areas aligned along $\mu$, see Figure~\ref{fig:closing_7x7x7_718_Bin_1.5}.  

\subsection{Enhancement of fault zones in compressed glass foam}
To characterise the structural behaviour of complex materials during loading, mechanical tests and simultaneous $\mu$CT imaging (in situ CT) can be used to estimate local displacements of the material on the micro-scale.
We analyse a glass foam which consists of very thin struts. 

\begin{figure}
    \centering
    \begin{subfigure}[b]{0.48\linewidth}
    \includegraphics[width=\linewidth]{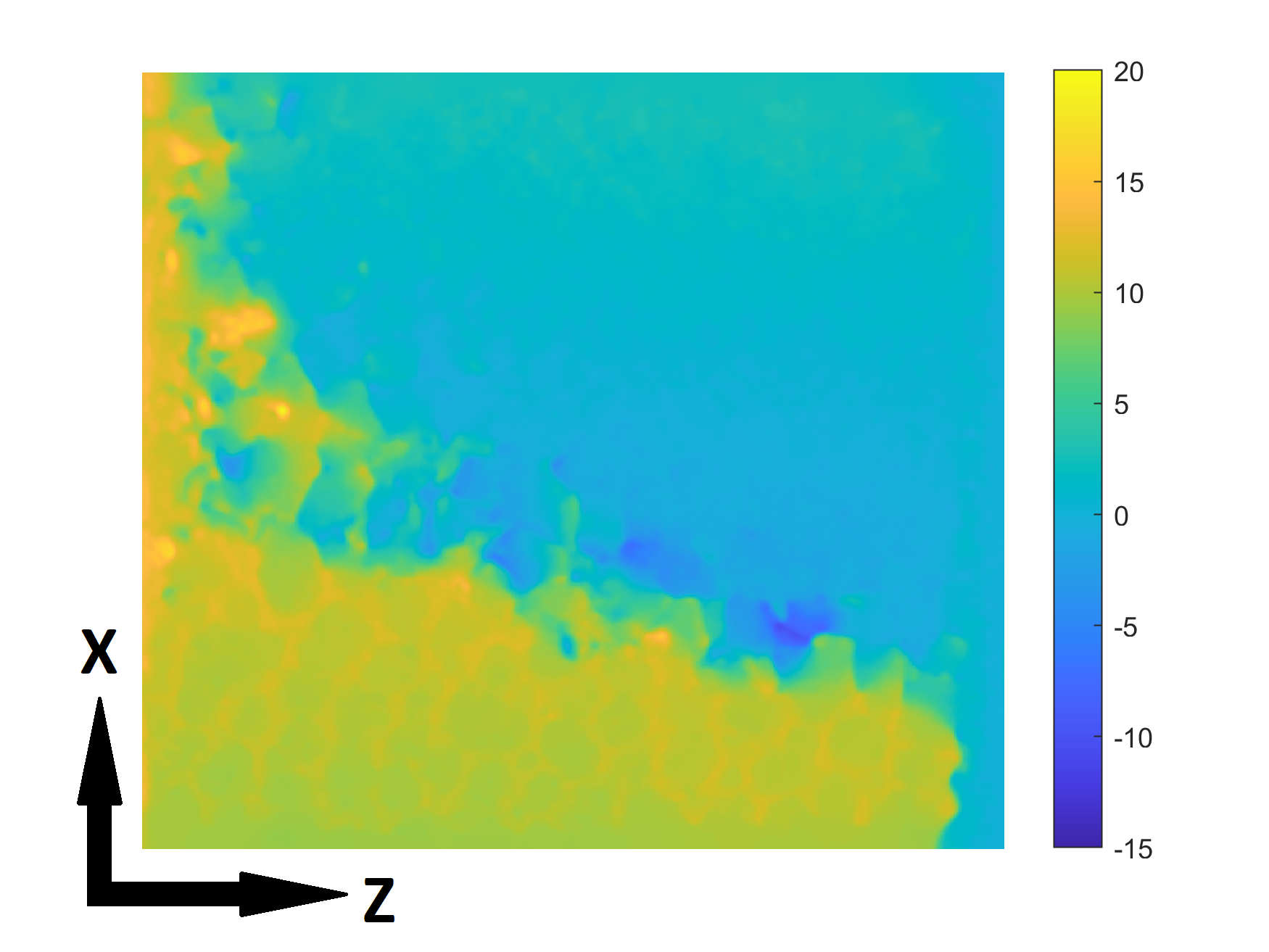}
    \caption{Displacement field}\label{fig:displacemnt_filed_w}
    \end{subfigure}
    \begin{subfigure}[b]{0.48\linewidth}
    \includegraphics[width=\linewidth]{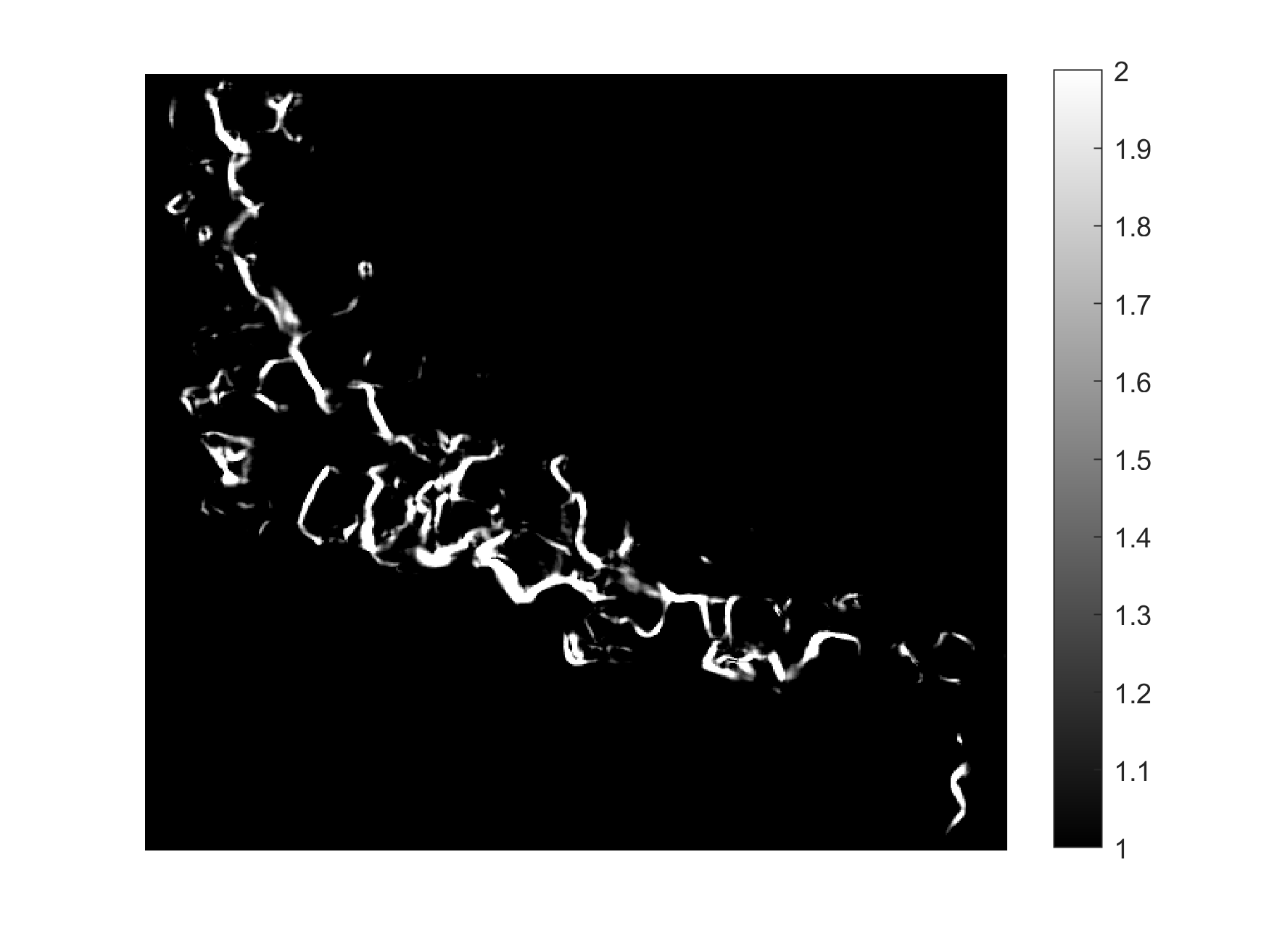}
    \caption{Morphological gradient}\label{fig:displacemnt_filed_w_gradient}
    \end{subfigure}
    
    \begin{subfigure}[b]{0.48\linewidth}
    \includegraphics[width=\linewidth]{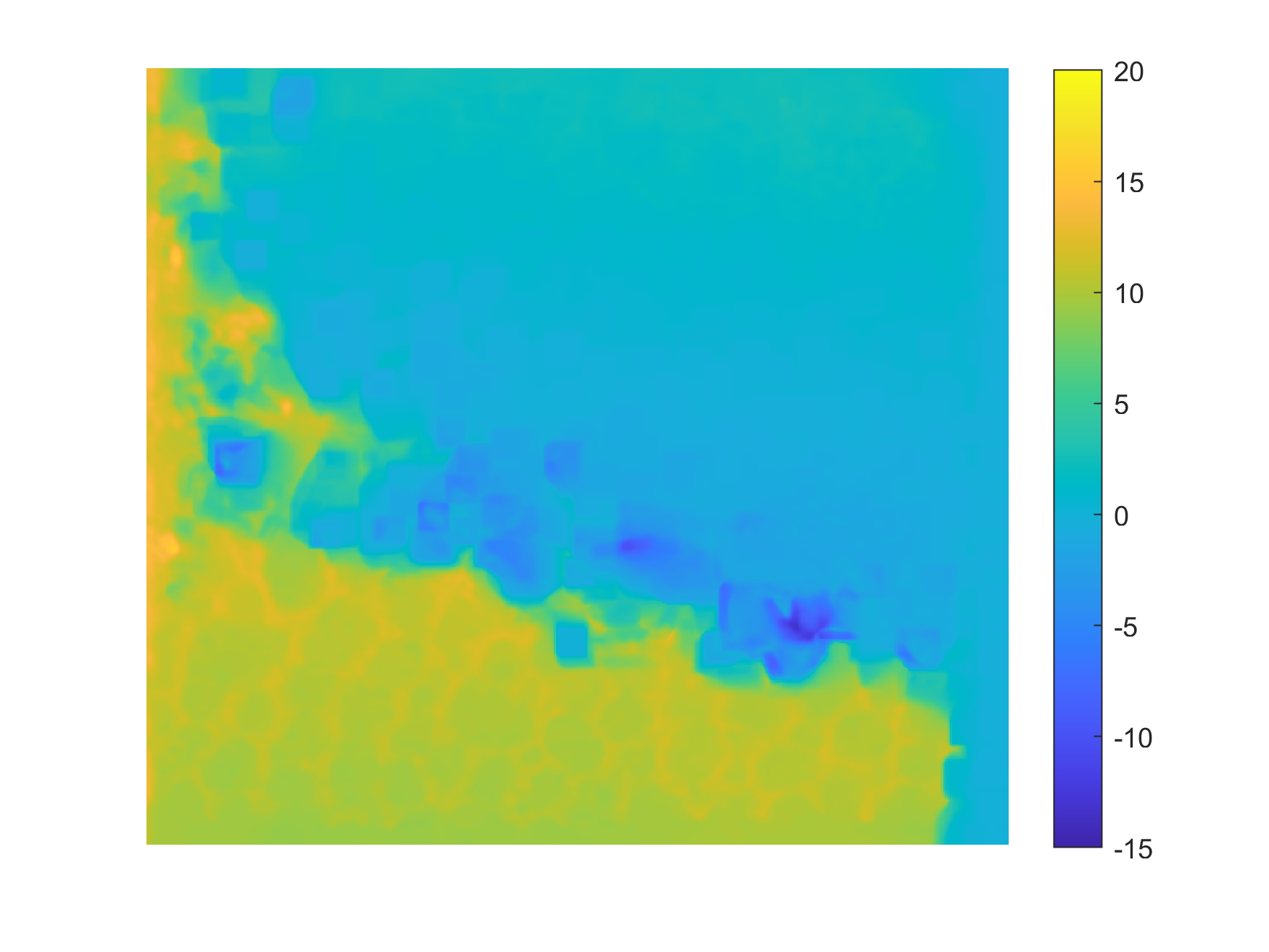}
    \caption{Erosion}\label{fig:displacemnt_filed_w_erosion}
    \end{subfigure}
    \begin{subfigure}[b]{0.48\linewidth}
    \includegraphics[width=\linewidth]{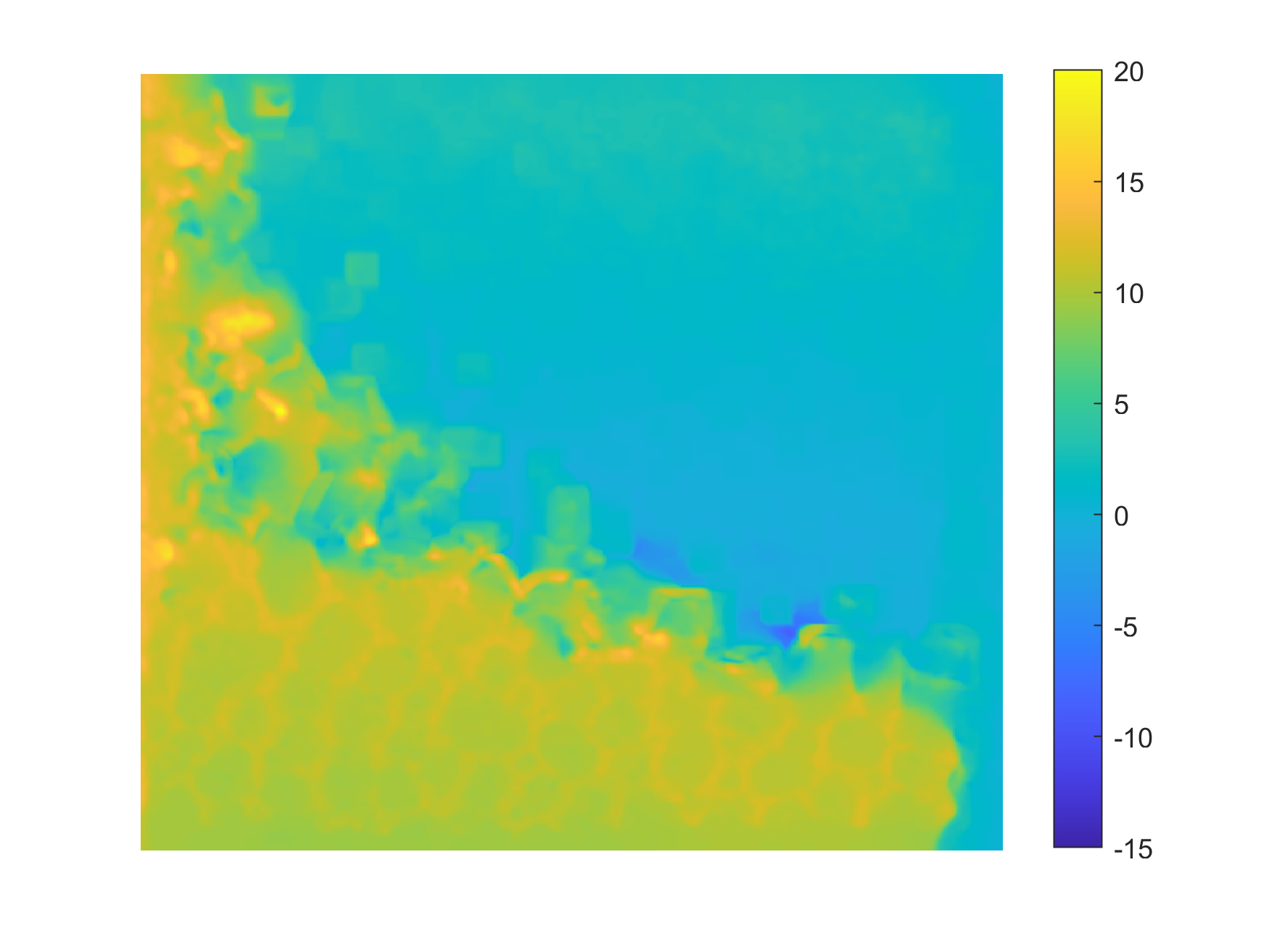}
    \caption{Dilation}\label{fig:displacemnt_filed_w_dilation}
    \end{subfigure}
    
    \begin{subfigure}[b]{0.48\linewidth}
    \includegraphics[width=\linewidth]{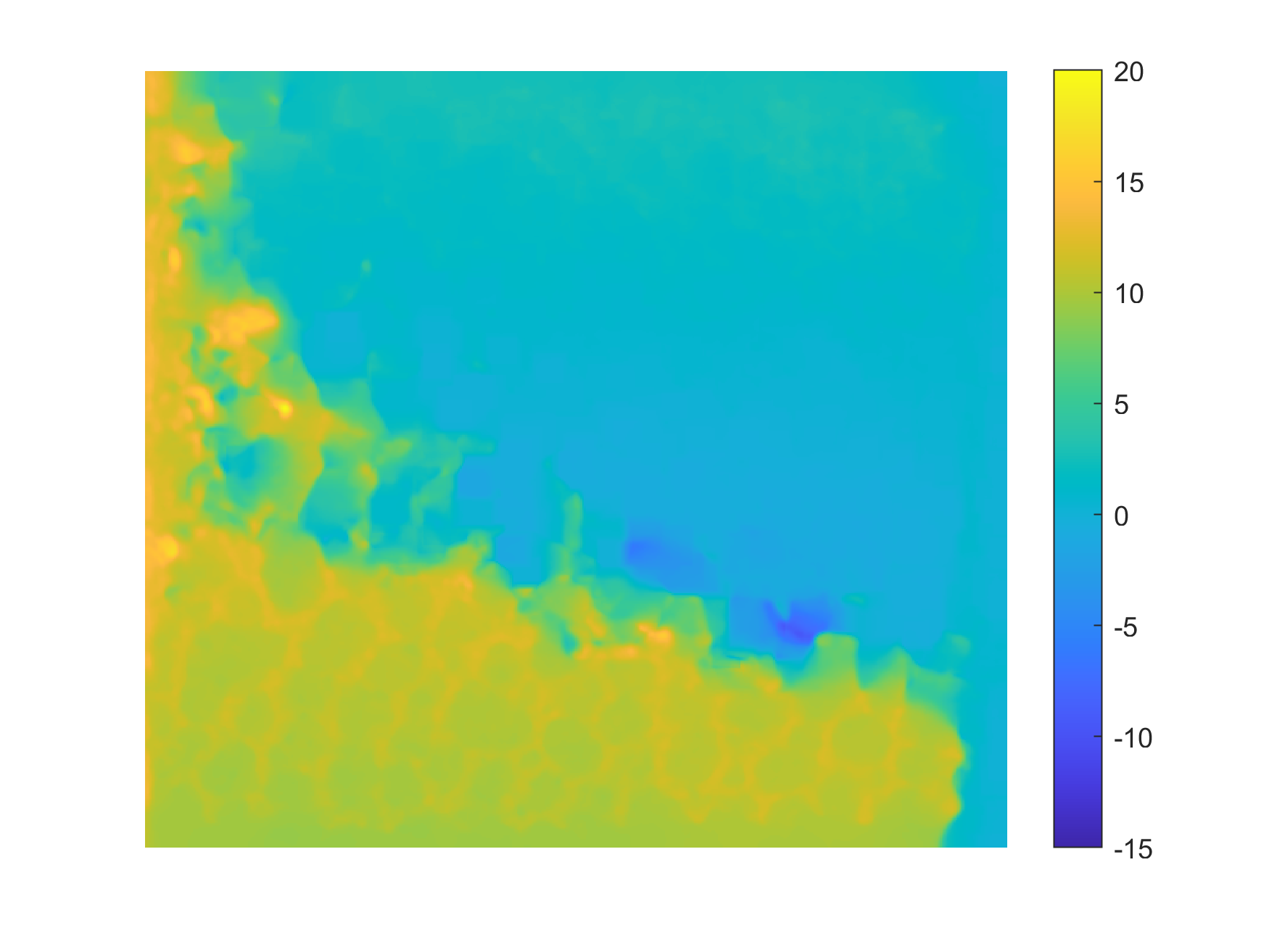}
    \caption{Opening}\label{fig:displacemnt_filed_w_closing}
    \end{subfigure}
    \begin{subfigure}[b]{0.48\linewidth}
    \includegraphics[width=\linewidth]{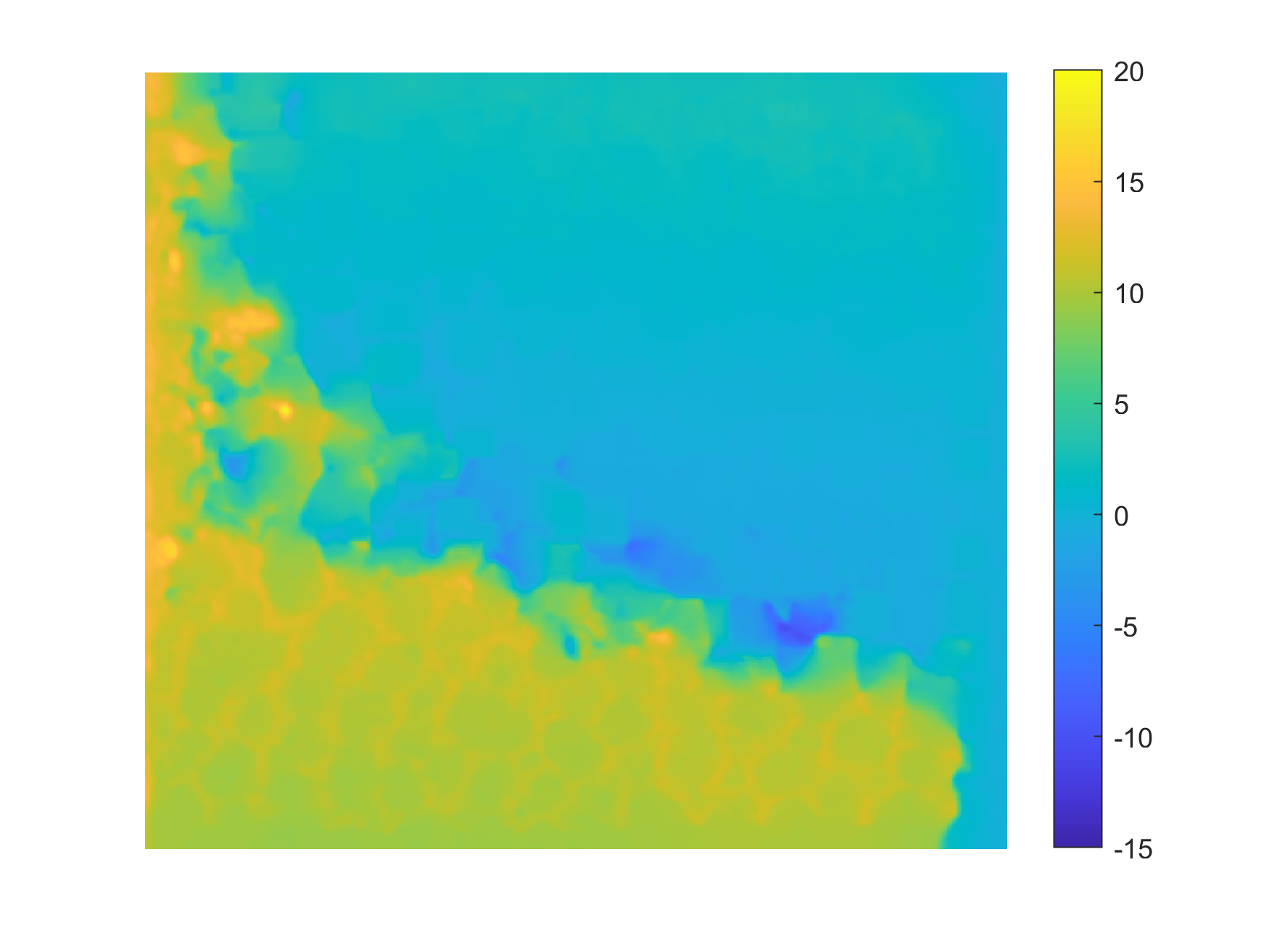}
    \caption{Closing}\label{fig:displacemnt_filed_w_opening}
    \end{subfigure}
    
    \caption{XZ-slices of displacement field based on \cite{Nogatz2021}. 
    Yellow colour indicates movement along the z-direction, and blue refers to the opposite direction. 
    The computed displacement field reflects the influence on the microstructure during compression (a).
    The morphological gradient (with $B$ a $3\times3$ square) enhances the fault zone (b). 
    The results of other morphological operators are given in (c-f).
    We chose $B$ as a $15\times15$ square.
    }
    \label{fig:Tessa_fig}
\end{figure}

The $\mu$CT images were recorded while the glass foam was compressed in the z-direction. 
The foam is expected to fail very suddenly due to its thin struts. 
Nogatz et al. \cite{Nogatz2021} computed the displacement field for the transition from strain level 1\% to 3.8\% which is shown in Figure \ref{fig:displacemnt_filed_w}.
During compression, a fault zone forms which is now post-processed by directional morphology. To this end, our operators are applied to the directional component of the displacement vectors. That is, the vectors were normalised for the calculation of the morphological operators and were subsequently scaled back to their original length. 
As we mostly expect compression along the loading direction, we choose $\mu=(0,0,1)$.

We enhance the fault zone with the morphological gradient, see Figure \ref{fig:displacemnt_filed_w_gradient}.
If motion estimation algorithms fail to reconstruct these sharp edges, we can enhance them by erosion, as shown in Figure \ref{fig:displacemnt_filed_w_erosion}. 
Some other materials however are known to show creep before failure. 
Here, a smoother transition seems more reasonable, which is obtained by a dilation, see Figure \ref{fig:displacemnt_filed_w_dilation}. 

Closing (Figure \ref{fig:displacemnt_filed_w_closing}) and opening (Figure \ref{fig:displacemnt_filed_w_opening}) remove misalignment in the transition between the two regions.
\section{Conclusion}\label{sec:Conclusion}
We have defined morphological operators for $\Sd$-valued images.
The required ordering for unit vectors is derived from the angular projection depth which enables a sound definition of $\Sd$-valued morphological operators and filters.
Relations of these morphological operators to their grey-scale counterparts are emphasised.
Additionally, we segmented misaligned fiber regions of glass fiber reinforced polymers and enhanced the fault region in a compressed glass foam.

Future work could address the removal of disturbances or the highlighting of directional changes in DTI images by morphological operations \cite{ZHANG2020-Overview,ZHANG2021-Alzheimer}.

\section*{Acknowledgment}
We thank Tessa Nogatz for providing the displacement field of the glass foam and Tin Bari\v{s}in for providing the glass fibre reinforced polymers images.

\section*{Appendix}
Given \eqref{eq:btsd(i) dot+ bssd(i)} the structuring element $\btsd$ fulfils the semi-group property \eqref{eq:structuring_function_semi-group} for 
\begin{align}
    \frac{\lVert i \rVert _2^2}{t},\frac{\lVert i \rVert _2^2}{s} < \pi      
    \label{eq:norm_t_s}
\end{align}
as follows:
Let $\overline{oi}$ be the line segment from the origin $o\in \Rq$ to $i\in \Rq$.
For $t,s > 0$, let $j^*$ be a point on $\overline{oi}\subset \Rq$ such that 
\begin{align}
	\lVert j^* \rVert_2 	&= \frac{s}{t+s} \lVert i \rVert _2 \label{eq:||j^*||_2} \\
	\lVert i-j^* \rVert_2 &= \frac{t}{t+s} \lVert i \rVert _2. \label{eq:||i-j^*||_2}	
\end{align}
Then,
\begin{align}
	(\alpha_t \Box \alpha_s)(i) 
	&\overset{\eqref{eq:infimal convolution}}{=} \inf_{j\in \Rq} \left\{ \alpha_t(i-j) + \alpha_s(j) \right\}\nonumber\\
    &\overset{\eqref{eq:alpha}}{=} \inf_{j\in \Rq} \left\{ \min\left(\frac{\lVert i-j \rVert _2^2}{t}, \pi \right) + \min\left(\frac{\lVert j \rVert _2^2}{s}, \pi \right) \right\}\nonumber\\
    &\overset{\eqref{eq:norm_t_s}}{=} \inf_{j\in \Rq} \left\{ \frac{||i-j||_2^2}{t} + \frac{\lVert j \rVert _2^2}{s} \right\}\nonumber\\
	&\overset{j^*\in \overline{oi}}{=} \frac{\lVert i-j^* \rVert_2^2}{t} + \frac{\lVert j^* \rVert_2^2}{s} \label{eq:j^*on-line}\\
	&\overset{\eqref{eq:||j^*||_2},\eqref{eq:||i-j^*||_2}}{=} \frac{t\lVert i \rVert _2^2}{(t+s)^2} + \frac{s\lVert i \rVert _2^2}{(t+s)^2}	\nonumber\\
	&= \frac{\lVert i \rVert _2^2}{(t+s)}\nonumber\\
    &\overset{\eqref{eq:norm_t_s}}{=}\min\left(\frac{\lVert i \rVert _2^2}{(t+s)}, \pi \right)\nonumber\\
	&\overset{\eqref{eq:alpha}}{=}\alpha_{t+s}(i). 
	\label{eq:angle infimal convolution}
\end{align}
Thus,
\begin{align*}
	\btsd(i) \dot{+} \bssd(i) 
	&\overset{\eqref{eq:btsd(i) dot+ bssd(i)}}{=} \Rmat((\alpha_t \Box \alpha_s)(i) ) 
        \overset{\eqref{eq:angle infimal convolution}}{=} \Rmat(\alpha_{t+s}(i) ) 
        \overset{\eqref{eq:bt-def}}{=}\btssd (i).
\end{align*}

\bibliographystyle{unsrt}  
\bibliography{references}

\end{document}